\documentclass[a4paper,10pt]{article}
\usepackage[english]{babel}
\usepackage[utf8]{inputenc}

\usepackage{amsmath}
\usepackage{amsfonts}
\usepackage{amssymb}
\usepackage{amsthm}
\usepackage{graphicx}
\usepackage{pdfpages}

\newcommand{\mytext}[1]{ \: \textrm{#1} \: }
\newcommand{\mysetdescr}[2]{\left\{ #1 \: \left| \: #2 \right. \right\} }

\newcommand{\myN}{\mathbb{N}}
\newcommand{\myNk}[1]{\underline #1}
\newcommand\mytimes{{\times}}

\newcommand\myurbild[1]{{#1}^{-1}}

\def\A{{\cal A}}
\def\E{{\cal E}}
\def\G{{\cal G}}
\def\H{{\cal H}}
\def\I{{\cal I}}
\def\J{{\cal J}}
\def\L{{\cal L}}
\def\S{{\cal S}}
\def\T{{\cal T}}

\newcommand{\strG}{\sqsubseteq_\Gamma}

\newcommand{\mygxiv}[2]{\Gamma_{#1}(#2)}
\newcommand{\gxiv}{\mygxiv{\xi}{v}}
\newcommand{\grxiv}{\mygxiv{\rho(\xi)}{v}}

\newcommand{\grgxiv}{\mygxiv{\rho_G(\xi)}{v}}

\newcommand{\agxi}{\alpha_{G,\xi}}

\newcommand{\myagxiv}[3]{\alpha_{#1, #2}(#3)}

\newcommand{\myagexiv}[3]{\alpha_{#1, \eta(#2)}(#3)}

\newcommand{\myaxiv}[2]{\alpha_{#1}(#2)}

\newcommand{\myaexiv}[2]{\alpha_{\eta(#1)}(#2)}

\newcommand{\agxiv}{\myagxiv{G}{\xi}{v}}
\newcommand{\agxiw}{\myagxiv{G}{\xi}{w}}

\newcommand{\axiv}{\myaxiv{\xi}{v}}
\newcommand{\axiw}{\myaxiv{\xi}{w}}

\newcommand{\myaiav}[2]{\alpha_{\iota({#1})}(#2)}

\newcommand{\eps}{\epsilon}

\newcommand{\xiv}{\xi(v)}
\newcommand{\xiw}{\xi(w)}

\newcommand{\myaTgxiv}[4]{\alpha^{#1}_{#2, #3}(#4)}
\newcommand{\myaTgxi}[3]{\alpha^{#1}_{#2, #3}}

\newcommand{\myaTxiv}[3]{\alpha^{#1}_{#2}(#3)}
\newcommand{\myaTxi}[2]{\alpha^{#1}_{#2}}

\newcommand{\aRgxiv}{\myaTgxiv{R}{G}{\xi}{v}}
\newcommand{\aRgxi}{\myaTgxi{R}{G}{\xi}}
\newcommand{\aRxiv}{\myaTxiv{R}{\xi}{v}}
\newcommand{\aRxi}{\myaTxi{R}{\xi}}
\newcommand{\aRgxiw}{\myaTgxiv{R}{G}{\xi}{w}}
\newcommand{\aRxiw}{\myaTxiv{R}{\xi}{w}}

\newcommand{\aSgexiv}{\myaTgxiv{S}{G}{\eta(\xi)}{v}}

\newcommand{\aSexiv}{\myaTxiv{S}{\eta(\xi)}{v}}
\newcommand{\aSexi}{\myaTxi{S}{\eta(\xi)}}

\newcommand{\myrhogxiv}[3]{\rho_{#1}(#2)(#3)}
\newcommand{\myrhogxi}[2]{\rho_{#1}(#2)}
\newcommand{\myrhoxiv}[2]{\rho(#1)(#2)}
\newcommand{\myrhoxi}[1]{\rho(#1)}

\newcommand{\rhogxiv}{\myrhogxiv{G}{\xi}{v}}
\newcommand{\rhogxi}{\myrhogxi{G}{\xi}}

\newcommand{\rhoxi}{\myrhoxi{\xi}}

\newcommand{\mf}[1]{\mathfrak{ #1 }}
\newcommand{\fa}{\mf{a}}
\newcommand{\fb}{\mf{b}}
\newcommand{\fc}{\mf{c}}
\newcommand{\fC}{\mf{C}}
\newcommand{\fd}{\mf{d}}
\newcommand{\fD}{\mf{D}}
\newcommand{\fP}{\mf{P}}
\newcommand{\fS}{\mf{S}}
\newcommand{\fT}{\mf{T}}
\newcommand{\fU}{\mf{U}}

\newcommand{\fTa}{\fT_a}

\newcommand{\Nin}{N^{in}}
\newcommand{\Nout}{N^{out}}

\newcommand{\rpha}{\rho(\phi)(\fa)}

\newcommand{\myfaan}[2]{\{ {#1}_1 \} \cup {#1}_{#2}}
\newcommand{\faaz}{\myfaan{\fa}{2}}

\newcommand{\fbbd}{\myfaan{\fb}{3}}

\newcommand{\myegxi}[2]{\eta_{#1}(#2)}
\newcommand{\exi}{\eta(\xi)}

\newcommand{\exiw}{\eta(\xi)(w)}

\newcommand{\sarie}{G \in \fD$, $\xi \in \S(G,R)$, $v \in V(G)}

\newcommand{\sarieStr}{G \in \fD'$, $\xi \in \S(G,R)$, $v \in V(G)}
\newcommand{\sariegxiStr}{G \in \fD'$, $\xi \in \S(G,R)}

\newcommand{\sarier}{G \in \fD$, $\xi \in \S(G,R)$, $v \in V(G)}

\newcommand{\arierStr}{G \in \fD'$, $\xi \in \H(G,R)$, $v \in V(G)}

\newcommand{\sarierStr}{G \in \fD'$, $\xi \in \S(G,R)$, $v \in V(G)}
\newcommand{\sariergxiStr}{G \in \fD'$, $\xi \in \S(G,R)}

\newcommand{\myXmn}[2]{X_{#1}^{#2}}
\newcommand{\myXmnC}[2]{\myXmn{\# #1}{\# #2}}
\newcommand{\Xmn}{\myXmn{m}{n}}

\def\BP{\begin{proof}}
\def\EP{\end{proof}}

\DeclareMathOperator{\id}{id}
\DeclareMathOperator{\Aut}{Aut}

\begin{document}

\theoremstyle{plain}
\newtheorem{condition}{Condition}
\newtheorem{theorem}{Theorem}
\newtheorem{definition}{Definition}
\newtheorem{corollary}{Corollary}
\newtheorem{lemma}{Lemma}
\newtheorem{proposition}{Proposition}

\title{\bf Generalized One-to-One Mappings between Homomorphism Sets of Digraphs}
\author{\sc Frank a Campo}
\date{\small Seilerwall 33, D 41747 Viersen, Germany\\
{\sf acampo.frank@gmail.com}}

\maketitle

\begin{abstract}
\noindent Structural properties of finite digraphs $R$ and $S$ are studied which enforce $\# \H(G,R) \leq \# \H(G,S)$ for every finite digraph $G \in \fD'$, where $\H(G,H)$ is the set of homomorphisms from $G$ to $H$, and $\fD'$ is a class of digraphs. In a previous study, we have seen that the key for such a relation between $R$ and $S$ is the existence of a {\em strong S-scheme} from $R$ to $S$. Such an S-scheme $\rho$ defines a one-to-one mapping $\rho_G : \S(G,R) \rightarrow \S(G,S)$ for every $G \in \fD'$, where $\S(G,H)$ is the set of homomorphisms from $G$ to $H$ mapping proper arcs of $G$ to proper arcs of $H$. In the present article, we characterize S-schemes $\rho$ which are induced by strict homomorphisms $\eps : \E(R) \rightarrow \E(S)$ between auxiliary systems of $R$ and $S$, and we analyze the mutual dependency between the properties of $\rho$ and $\eps$. Wide applicability of the theory is ensured by specifying the auxiliary systems $\E(R)$ and $\E(S)$ as {\em EV-systems} of $R$ and $S$. The results are applied on a rearrangement method for digraphs and on undirected graphs.
\newline

\end{abstract}

\section{Introduction} \label{sec_introduction}

The number of homomorphisms between directed graphs (digraphs) may carry important information about structure. Freedman et al.\ \cite{Freedman_etal_2007} characterized in 2007 graph parameters which can be expressed as numbers of homomorphisms into weighted graphs. The still open reconstruction conjecture asks in different fields of graph theory \cite{Hell_Nesetril_2004,Schroeder_2016}, if two objects with at least four vertices are isomorphic if all numbers of embeddings of certain subgraphs into them are equal. An early classical result is the Theorem of Lov\'{a}sz \cite{Lovasz_1967} from 1967 which states that numbers of homomorphisms distinguish non-isomorphic ``relational structures''. With $\H(G,H)$ defined as the set of homomorphisms from a digraph $G$ to a digraph $H$, the following specifications are relevant for our purpose:

\begin{theorem}[Lov\'{a}sz \cite{Lovasz_1967}] \label{theo_Lovasz_original}
Let $\fC$ be the class of finite digraphs or the class of finite posets. Then, for $R, S \in \fC$
\begin{align*}
R & \simeq S \\
\Leftrightarrow \; \; \; \# \H(G,R) & = \# \H(G,S) \; \mytext{for every } G \in \fC.
\end{align*}
For the class of finite posets, the equivalence holds also if we replace the homomorphism sets by the sets of strict order homomorphisms.
\end{theorem}

For digraphs, a short and simple proof of the theorem is contained in \cite{Hell_Nesetril_2004} which - with minor modification - works for posets, too.

The infinite vector $\L(H) \equiv ( \# \H(G,H) )_{G \in \fC}$ is called the {\em Lov\'{a}sz-vector of $H$}.  In the last decades, topics related to it have found interest \cite{Borgs_etal_2006,Lovasz_2006,Freedman_etal_2007,
Borgs_etal_2008,Lovasz_Szegedy_2008,Schrijver_2009,
Borgs_etal_2012,Cai_Govorov_2020} in connection with vertex and edge weights. For undirected graphs, Dvo\v{r}\'{a}k \cite{Dvorak_2010} investigated in 2010 proper sub-classes $\fU'$ of undirected graphs for which the partial Lov\'{a}sz-vector $( \# \H(G,H) )_{G \in \fU'}$ still distinguishes graphs; the distinguishing power of the vector $( \# \H(G,H) )_{H \in \fU'}$ has been investigated by Fisk \cite{Fisk_1995} in 1995.

This paper continues the work of the author about the pointwise less-equal-relation between partial Lov\'{a}sz-vectors of digraphs:
\begin{quote}
{\em Given a class $\fD'$ of digraphs, what is it in the structure of digraphs $R$ and $S$ that enforces}
\end{quote}
\begin{equation} \label{fragestellung}
\# \H(G,R) \leq \# \H(G,S) \; \mytext{\em for every } G \in \fD' \mytext{\em ?}
\end{equation}

\begin{figure} 
\begin{center}
\includegraphics[trim = 75 710 185 70, clip]{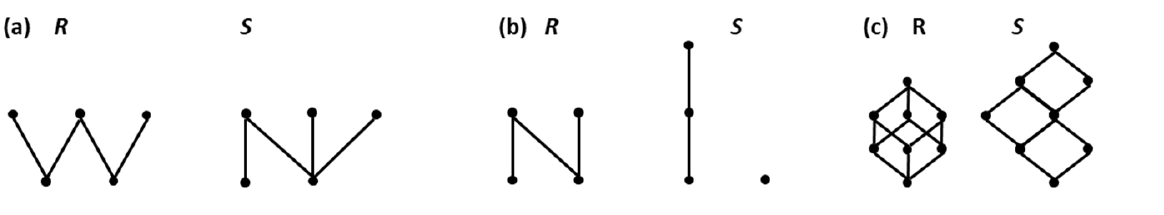}
\caption{\label{figure_Intro} Hasse-diagrams of three pairs $R$, $S$ of posets with $\# \H(P,R) \leq \# \H(P,S)$ for every finite poset $P$.}
\end{center}
\end{figure}

The starting point of the work was the pair of posets $R$ and $S$ in Figure \ref{figure_Intro}(c). The author \cite[Theorem 5]{aCampo_2018} has proven that for these posets we have $\# \H(P,R) \leq \# \H(P,S)$ for every finite poset $P$. Additional non-trivial examples for the relation ``$\# \H(P,R) \leq \# \H(P,S)$ for every finite poset $P$'' are shown in the Figures \ref{figure_Intro}(a)-(b); more pairs of digraphs fulfilling \eqref{fragestellung} are contained in the Figures \ref{figure_RS_nichtRearr} and \ref{fig_TableConstr}.

In the recent paper, a first attempt is made to tackle the theoretical aspect of the question. After the preparatory Section \ref{sec_preparation}, a main result from \cite{aCampo_toappear_0} is recalled in Section \ref{subsec_DefHomSchemes}: it has been shown \cite[Theorem 2]{aCampo_toappear_0}, that for several important classes of digraphs $\fD'$ and $R \in \fD'$, \eqref{fragestellung} is implied by
\begin{align} \label{fragestellung_strict}
\# \S(G,R) \; \leq \; \# \S(G,S) \; \; \mytext{for all} \; G \in \fD'
\end{align} 
where $\S(G,H)$ is the set of {\em strict} homomorphisms from $G$ to $H$, i.e., of those homomorphisms from $G$ to $H$ mapping proper arcs of $G$ to proper arcs of $H$.

Due to the implication $\eqref{fragestellung_strict} \Rightarrow \eqref{fragestellung}$, our interest focuses on the existence of a generalized mapping $\rho$ providing a mapping $\rho_G : \S(G,R) \rightarrow \S(G,S)$ for every $G \in \fD'$. (A formally satisfying definition is given in Definition \ref{def_S_scheme}.) We call such a generalized mapping an {\em S-scheme from $R$ to $S$}, and we call $\rho$ {\em strong} iff $\rho_G$ is one-to-one for every $G \in \fD'$.

For a (strong) S-scheme $\rho$ and $G, G' \in \fD'$ with $G \not= G'$, the mappings $\rho_G$ and $\rho_{G'}$ are in principle independent from each other, no matter how closely $G$ and $G'$ are related. Even in the case of $\xi \in \S(G,R) \cap \S(G',R)$, the resulting homomorphisms $\rho_G(\xi) \in \S(G,S)$ and $\rho_{G'}(\xi) \in \S(G',S)$ must not be related to each other in any way. It are the {\em induced S-schemes} defined in Section \ref{subsec_simple_induced} which introduce some regularity by ``mapping similar things to similar things'', and they are the subject of this paper. In particular, S-schemes are suitable for constructive approaches, and we will pay special attention to the types of regularity introduced by them.

For an induced S-scheme $\rho$ from $R$ to $S$, we have for all $G \in \fD'$
\begin{equation*}
\rho_G \; = \; \phi \circ \eps \circ \alpha_G,
\end{equation*}
where $\alpha$ is a simple S-scheme from $R$ to an auxiliary digraph $\E(R)$, $\phi$ is a strict homomorphism from another auxiliary digraph $\E(S)$ to $S$, and $\eps : \E(R) \rightarrow \E(S)$ is a strict homomorphism. Obviously, the properties of $\rho$ are mainly determined by the digraphs $\E(R)$ and $\E(S)$ and the homomorphism $\eps$ between them. From a theoretical point of view it is interesting to analyze the mutual dependency between the properties of $\rho$ and $\eps$, and from a practical point of view, there is a perspective to construct a strong S-scheme $\rho$ from $R$ to $S$ with desired properties by designing (hopefully simple) digraphs $\E(R)$ and $\E(S)$ and an appropriate homomorphism $\eps$ between them.

The theoretical aspect is in the focus in what follows. Still in Section \ref{subsec_simple_induced}, we characterize the induced S-schemes in Theorem \ref{theo_induced_exists}. In Theorem \ref{theo_eta_RS}(1), we show that a regularity condition on $\eps$ called ``Condition \ref{cond_univ_aexiv}'' is equivalent to a certain regular behavior of the induced S-scheme $\rho$, and in Theorem \ref{theo_eta_RS}(2), we see that an additional property of $\eps$ forces $\rho$ to be strong.

In Section \ref{subsec_EVSystems_Props}, we specify $\E(R)$ and $\E(S)$ as {\em EV-systems} and describe their dependency on the class of digraphs $\fD'$ they are referring to. In Theorem \ref{theo_eta_RS_inv} in Section \ref{subsec_inv_Theo_Eta_RS_2}, we use EV-systems in showing that for the classes of digraphs we are mainly interested in, Theorem \ref{theo_eta_RS}(2) can be inverted: if an induced strong S-scheme $\rho$ behaves sufficiently regular, then the corresponding $\eps$ is one-to-one and fulfills Condition \ref{cond_univ_aexiv}. Because Condition \ref{cond_univ_aexiv} is unwieldy to check, we show in Proposition \ref{prop_repl_Cond1} in Section \ref{subsec_replacement_Cond1} how it can be replaced by a more handy condition.
 
In Section \ref{subsec_rearr_induced}, we take up the rearrangement method developed in \cite{aCampo_toappear_0}. By means of this method, a digraph $R$ fulfilling certain conditions is rearranged in such a way that there exists a strong S-scheme $\rho$ from $R$ to the digraph $S$ resulting from the rearrangement. We see that $\rho$ is in fact an induced strong S-scheme, and we describe the corresponding homomorphism $\eps : \E(R) \rightarrow \E(S)$. In Section \ref{subsec_example}, we discuss in detail the posets in Figure \ref{figure_Intro}(a) and Figure \ref{figure_Intro}(b) under these view points.

Finally, in Section \ref{sec_undirected}, we transfer our concepts and results to undirected graphs.

\section{Basics and Notation} \label{sec_preparation}

A {\em (finite) directed graph} or {\em digraph} $G$ is an ordered pair $(V(G),A(G))$ in which $V(G)$ is a non-empty, finite set and $A(G) \subseteq V(G) \mytimes V(G)$ is a binary relation on $V(G)$. We write $vw$ for an ordered pair $(v,w) \in V(G) \mytimes V(G)$. The elements of $V(G)$ are called the {\em vertices} of $G$ and the elements of $A(G)$ are called the {\em arcs} of $G$. A digraph $G$ is {\em reflexive},  or {\em symmetric}, or {\em antisymmetric}, etc., iff the relation $A(G)$ has the respective property. A {\em partially ordered set (poset)} is a reflexive, antisymmetric, transitive digraph. 

For a digraph $G$ and a non-empty set $X \subseteq V(G)$, the {\em digraph $G \vert_X$ induced on $X$} is $( X, A(G) \cap ( X \mytimes X) )$. The {\em direct sum} $G + H$ of digraphs with disjoint vertex sets is defined as usual.

Vertices $v, w \in V(G)$ are {\em adjacent} iff $vw \in A(G)$ or $wv \in A(G)$. The {\em open neighborhood} $N_G(v)$ of $v \in V(G)$ is the set of all $w \in V(G) \setminus \{ v \}$ adjacent to $v$. Furthermore,
\begin{align*}
\Nin_G(v) & \; \equiv \; \mysetdescr{ w \in N_G(v) }{ wv \in A(G) }, \\
\Nout_G(v) & \; \equiv \; \mysetdescr{ w \in N_G(v) }{ vw \in A(G) }.
\end{align*}
An arc $vw \in A(G)$ is called {\em proper} iff $v \not= w$; otherwise, it is called a {\em loop}. All possible loops of $G$ are collected in the {\em diagonal} $\Delta_G \equiv \mysetdescr{(v,v)}{v \in V(G)}$. $G^* \equiv (V(G), A(G) \setminus \Delta_G)$ is the digraph $G$ {\em with loops removed}.

A sequence $v_0, \ldots , v_I$ of vertices of $G$ with $I \in \myN$ is called a {\em walk} iff $v_{i-1}v_i \in A(G)$ for all $1 \leq i \leq I$. The walk is {\em closed} iff $v_0 = v_I$. A digraph is {\em acyclic} iff it does not contain a closed walk.

Let $G$ be a digraph. With $\T$ denoting the set of all transitive relations $T \subseteq V(G) \mytimes V(G)$ with $A(G) \subseteq T$, the {\em transitive hull} $( V(G), \cap \T )$ of $G$ is the digraph with vertex set $V(G)$ and the (set-theoretically) smallest transitive arc set containing $A(G)$. 

Given digraphs $G$ and $H$, we call a mapping $\xi : V(G) \rightarrow V(H)$ a {\em homomorphism} from $G$ to $H$ if $\xi(v) \xi(w) \in A(H)$ for all $vw \in A(G)$. For such a mapping, we write $\xi : G \rightarrow H$; we collect the homomorphisms in the set
\begin{align*}
\H(G,H) & \; \equiv \; \mysetdescr{ \xi : V(G) \rightarrow V(H) }{ \xi \mytext{ is a homomorphism} }.
\end{align*}
$\Aut(G)$ is the set of automorphisms of a digraph $G$. Isomorphism is indicated by ``$\simeq$''.

Every homomorphism $\xi : G \rightarrow H$ maps loops in $G$ to loops in $H$, but proper arcs of $G$ can be mapped to both, loops and proper arcs of $H$. We call a homomorphism from $G$ to $H$ {\em strict} iff it maps all proper arcs of $G$ to proper arcs of $H$.
\begin{align*}
\S(G,H) & \; \equiv \; \H(G,H) \cap \H(G^*,H^*)
\end{align*}
is the set of strict homomorphisms from $G$ to $H$. The set $\H(G^*,H^*) \setminus \H(G,H)$ contains all homomorphisms from $G^*$ to $H^*$ which map a vertex belonging to a loop in $G$ to a vertex of $H$ not belonging to a loop. A mapping $\xi : V(G) \rightarrow V(H)$ is thus a strict homomorphism, iff it maps loops in $G$ to loops in $H$ and proper arcs of $G$ to proper arcs of $H$. If $H$ is reflexive, then $\S(G,H) = \H(G^*,H^*)$; for posets $P$ and $Q$, the set $\S(P,Q) = \H(P^*,Q^*)$ is the set of strict order homomorphisms from $P$ to $Q$.

In order to avoid unnecessary formalism, we work with representative systems of classes of digraphs with respect to isomorphism and not with the classes itself. $\fD$ is a representative system of the class of all digraphs with finite non-empty vertex set. Furthermore,
\begin{align*}
\fTa & \; \equiv \; \mysetdescr{ G \in \fD }{ G^* \mytext{is acyclic} }, \\
\fP \; & \; \equiv \; \mysetdescr{ P \in \fTa }{ P \mytext{is a poset} }, \\
\fP^* & \; \equiv \; \mysetdescr{ P^* }{ P \in \fP }.
\end{align*}
Equivalently, $\fTa$ can be characterized as an representative system of the class of digraphs with antisymmetric transitive hull (which is the reason for the choice of the symbol $\fTa$) or as the class of subgraphs of posets. $\fP^*$ is the class of posets with loops removed, i.e., the class of finite irreflexive antisymmetric transitive digraphs. $\fP^*$ is of interest for us, because every result about homomorphism sets $\H(P,Q)$ with $P, Q \in \fP^*$ directly translates into a result about the sets of strict order homomorphisms between posets and vice versa. For example, the addendum in Theorem \ref{theo_Lovasz_original} says that the stated equivalence also holds for $\fC = \fP^*$.

We assume $\fP^* \subset \fT_a$, and we will tacitly assume that every digraph we construct in what follows belongs to the respective representative system without exchange of vertices. Nevertheless, for $G, H \in \fD$, we retain the notation $G \simeq H$ instead of $G = H$ in order to emphasize that it is structural equivalence we are dealing with, not physical identity.

From set theory, we use additionally the following notation:
\begin{align*}
\myNk{0} & \equiv  \emptyset, \\
\myNk{n} & \equiv  \{ 1, \ldots, n \} \mytext{for every} n \in \myN.
\end{align*}

$\A(X,Y)$ is the set of mappings from $X$ to $Y$. For $f \in \A(X,Y)$ and $X' \subseteq X$, we write $f \vert_{X'}$ for the pre-restriction of $f$ to $X'$, and for $Y' \subseteq Y$ with $f(X) \subseteq Y'$ we write $f \vert^{Y'}$ for the post-restriction of $f$ to $Y'$. Furthermore, we use the symbol $\myurbild{f}(Y'')$ for the pre-image of $Y'' \subseteq Y$ under $f$; for $y \in Y$, we simply write $\myurbild{f}(y)$ instead of $\myurbild{f}( \{y\} )$. However, in the proof of Lemma \ref{lemma_eps_nicht_1t1_Cond1}, we use the symbol $\beta^{-1}$ also for the inverse of a bijective mapping $\beta$. $\id_X$ is the identity mapping of a set $X$.

Finally, we use the {\em Cartesian product}. Let $\I$ be a non-empty set, and let $M_i$ be a non-empty set for every $i \in \I$. Then the Cartesian product of the sets $M_i, i \in \I$, is defined as
\begin{eqnarray*}
\prod_{i \in \I} M_i & \; \equiv \; & 
\mysetdescr{ f \in \A \big( \I, \bigcup_{i \in \I} M_i \big)}{ f(i) \in M_i \mytext{for all} i \in \I }.
\end{eqnarray*}

\section{S-schemes} \label{sec_HomSchemes}

In Section \ref{subsec_DefHomSchemes}, we recall the main concepts and results from \cite{aCampo_toappear_0}. It turns out that for our purpose, so-called {\em S-schemes} are in the focus: generalized mappings $\rho$ providing a mapping $\rho_G : \S(G,R) \rightarrow \S(G,S)$ for every $G \in \fD' \subseteq \fD$. In Section \ref{subsec_simple_induced}, we introduce {\em induced S-schemes}. These S-schemes can be described effectively by two auxiliary digraphs and a strict homomorphism $\eps$ between them. We characterize the induced S-schemes and show how regularity conditions on $\eps$ translate into regular behavior of the induced S-scheme $\rho$ and vice versa.

\subsection{Recapitulation} \label{subsec_DefHomSchemes}

Let $R, S \in \fD$ and $\fD' \subseteq \fD$. Assume that there exists a one-to-one homomorphism $\sigma : R \rightarrow S$. Then, for every $G \in \fD'$ with $\H(G,R) \not= \emptyset$, we get a one-to-one mapping $r_G : \H(G,R) \rightarrow \H(G,S)$ by setting for every $\xi \in \H(G,R)$
\begin{align} \label{eq_rho_sigma}
r_G(\xi) & \equiv \sigma \circ \xi.
\end{align}
A one-to-one homomorphisms from $R$ to $S$ delivers thus a ``natural'' (or: trivial) example for $ \# \H(G,R) \leq \# \H(G,S)$ for every $G \in \fD'$. For the general investigation of this relation between $R$ and $S$, the author \cite{aCampo_2018,aCampo_toappear_0} has introduced the following concepts:

\begin{definition}[\cite{aCampo_toappear_0}, Definition 3] \label{def_Hom_scheme}
Let $\fD' \subseteq \fD$ be a representative system of a class of digraphs. For $R, S \in \fD$, we call a mapping
\begin{align*}
\rho & \; \in \prod_{G \in \fD'} \A( \H(G,R), \H(G,S) )
\end{align*}
a {\em Hom-scheme from $R$ to $S$ (with respect to $\fD'$)}, and we call it {\em strong} iff $\rho_G : \H(G,R) \rightarrow \H(G,S)$ is one-to-one for every $G \in \fD'$. We say that a Hom-scheme $\rho$ from $R$ to $S$ is a {\em $\Gamma$-scheme}, iff
\begin{align} \label{def_grxiv_gxiv}
\gxiv & = \grgxiv
\end{align}
holds for every $G \in \fD', \xi \in \H(G,R), v \in V(G)$, where $\gxiv$ is the connectivity component of $v$ in $\myurbild{\xi}(\xiv)$ and $\grgxiv$ is the connectivity component of $v$ in $\myurbild{\rhogxi}(\rhogxiv)$. We write $R \strG S$ iff a strong $\Gamma$-scheme from $R$ to $S$ exists. If $G$ is fixed, we write $\rho(\xi)$ instead of $\rho_G(\xi)$.
\end{definition}

Here as in the following, it does not matter if there is a $G \in \fD'$ with $\H(G,R) = \emptyset$; in this case, $\rho_G = ( \emptyset, \emptyset, \H(G,S))$. The (trivial) Hom-scheme \eqref{eq_rho_sigma} is always a strong $\Gamma$-scheme. We will generalize the concept of such simple Hom-schemes in Section \ref{subsec_simple_induced}.

A Hom-scheme $\rho$ from $R$ to $S$ is strong if we can determine $\xiv$ by means of $\rho(\xi)$ for all $\arierStr$. We say that we can {\em reconstruct $\xi$ by means of $\rho(\xi)$}. Obviously, we can reconstruct $\xi$ by means of $\rho(\xi)$ if we can determine $\myurbild{\xi}(v)$ for every $v \in R$ by means of $\rho(\xi)$.

A (strong) $\Gamma$-scheme is a (strong) Hom-scheme obeying the regularity condition \eqref{def_grxiv_gxiv} in mapping $\H(G,R)$ to $\H(G,S)$ for every $G \in \fD'$. This condition has been introduced in \cite{aCampo_toappear_0} under an application-oriented point of view. It may look like an additional difficulty posed upon a question difficult enough in itself. However, in fact, \eqref{def_grxiv_gxiv} is a regularity condition making things manageable by introducing structure. For a (strong) Hom-scheme $\rho$  from $R$ to $S$ and $G, G' \in \fD'$ with $G \not= G'$, the mappings $\rho_G$ and $\rho_{G'}$ are independent; even in the case of $\xi \in \H(G,R) \cap \H(G',R)$, there must be no similarity between the image-homomorphisms $\rho_G(\xi)$ and $\rho_{G'}(\xi)$. It are the $\Gamma$-schemes and in particular the induced S-schemes defined in Section \ref{subsec_simple_induced} which ensure that ``similar things are mapped to similar things''. In this way, S-schemes are suitable for constructive approaches, and we will pay particular attention to the question which type of regularity is introduced by them.

The regularity condition \eqref{def_grxiv_gxiv} is in particular plausible if we regard a Hom-scheme as a technical apparatus which assigns to every $\xi \in \H(G,R)$ a well-fitting $\rho(\xi) \in \H(G,S)$. If we allow $\gxiv \subset \grxiv$ for $x \in P$, then \cite[Lemma 1]{aCampo_toappear_0} tells us that $\rho(\xi)$ preserves the structure of $G$ around $v$ worse than $\xi$, which is not satisfying. And in the case $\gxiv \not\subseteq \grxiv$, $\rho(\xi)$ has to re-distribute the points of $\gxiv \setminus \grxiv \subseteq \gxiv \setminus \{ v \}$ in $S$. Because the sets $\gxiv \setminus \rhogxiv$ can be arbitrarily complicated, this re-distribution process may require many single case decisions, which is out of the scope of a technical apparatus.

It is easily seen \cite[Corollary 3]{aCampo_toappear_0} that a homomorphism $\xi \in \H(G,H)$ is strict iff $\gxiv = \{ v \}$ for all $v \in G$. Because a $\Gamma$-scheme preserves the sets $\gxiv$, we have for every Hom-scheme $\rho$ from $R$ to $S$
\begin{equation} \label{GScheme_strict}
\rho \mytext{ $\Gamma$-scheme } \quad \Rightarrow \quad 
\rho_G( \S(G,R) ) \subseteq \S(G,S) \quad \mytext{for all } G \in \fD',
\end{equation}
hence $\# \S(G,R) \leq \# \S(G,S)$ for all $G \in \fD'$ if $\rho$ is additionally strong. One of the main results of \cite{aCampo_toappear_0} adds the direction ``$\Leftarrow$'' to this implication:
\begin{theorem}[\cite{aCampo_toappear_0}, Theorem 2] \label{theo_GschemeOnStrict}
Let \begin{align*}
R & \in \fD' = \fD, \\
R & \in \fTa \subseteq \fD' \subseteq \fD, \\
\mytext{or} \quad R & \in \fD' \mytext{ with } \fD' = \fP \mytext{ or } \fD' = \fP^*.
\end{align*}
Then, for all $S \in \fD$, there exists a strong $\Gamma$-scheme from $R$ to $S$ with respect to $\fD'$, iff
\begin{equation*}
\# \S(G,R) \; \leq \; \# \S(G,S) \; \; \mytext{for all} \; G \in \fD'.
\end{equation*}
\end{theorem}

In the investigation of strong $\Gamma$-schemes, it is thus obvious to give special attention to the restriction of Hom-schemes to sets of strict homomorphisms:
\begin{definition} \label{def_S_scheme}
Let $R, S \in \fD$, $\fD' \subseteq \fD$. We call a mapping
\begin{align*}
\rho & \; \in \; \prod_{G \in \fD'} \A( \S(G,R), \S(G,S) )
\end{align*}
an {\em S-scheme from $R$ to $S$ with respect to $\fD'$}. We call an S-scheme from $R$ to $S$ {\em strong} iff the mapping $\rho_G : \S(G,R) \rightarrow \S(G,S)$ is one-to-one for every $\sariergxiStr$.
\end{definition}

In fact, we have seen:

\begin{proposition}[\cite{aCampo_toappear_0}, Proposition 2] \label{prop_extend_Sscheme}
In the constellations of $R$ and $\fD'$ described in Theorem \ref{theo_GschemeOnStrict}, a strong S-scheme $\rho$ from $R$ to $S$ can always be extended to a strong $\Gamma$-scheme $\rho'$ with $\rho'_G \vert_{\S(G,R)}^{\S(G,S)} = \rho_G$ for all $G \in \fD'$. 
\end{proposition}

Theorem \ref{theo_GschemeOnStrict} and Proposition \ref{prop_extend_Sscheme} remain valid also for other sub-classes $\fD'$ of $\fD$, e.g., for the digraphs (posets) with at most $k$ vertices or at most $k$ arcs, and for the class of digraphs in $\fTa$ for which the maximal length of a walk without loops is at most $k$.

\subsection{Induced S-schemes} \label{subsec_simple_induced}

\begin{figure} 
\begin{center}
\includegraphics[trim = 70 620 290 70, clip]{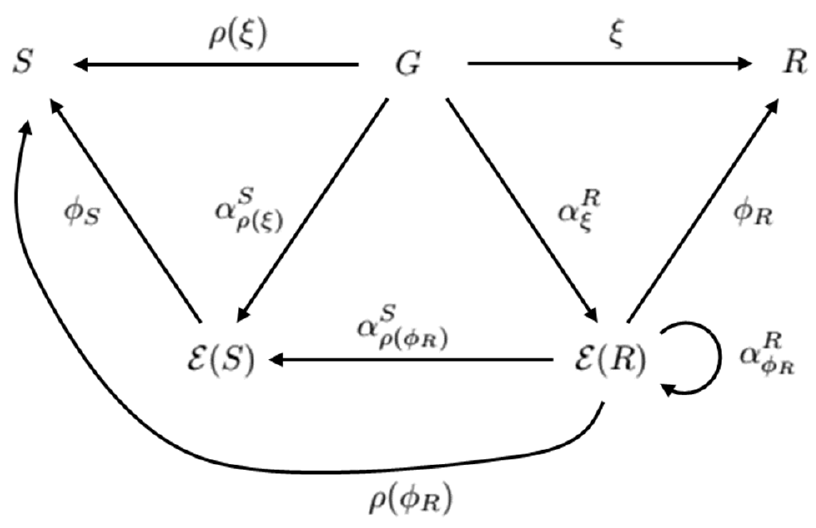}
\caption{\label{figure_KommDiagr} Homomorphisms and their images for an induced S-scheme $\rho$. Explanations in text.}
\end{center}
\end{figure}

In the following definition, the main property of the trivial strong $\Gamma$-scheme in \eqref{eq_rho_sigma} is generalized and transferred to S-schemes:

\begin{definition} \label{def_EVSys}
Let $T, \E(T) \in \fD$, and let $\alpha$ be an S-scheme from $T$ to $\E(T)$ with respect to $\fD' \subseteq \fD$. We call $\alpha$ a {\em simple S-scheme} iff there exists a strict homomorphism $\phi_T : \E(T) \rightarrow T$ with
\begin{equation} \label{defeq_EVsyst}
\forall \; G \in \fD', \xi \in \S(G,T) \mytext{:} \; \phi_T \circ \agxi \; = \; \xi.
\end{equation}
If $G$ is fixed, we write $\axiv$ instead of $\agxiv$. 
\end{definition}
Because $\alpha$ is an S-scheme, $\agxi$ is strict for all $\G \in \fD'$, $\xi \in \S(G,T)$. \eqref{defeq_EVsyst} shows that $\alpha_G : \S(G,T) \rightarrow \S(G,\E(T))$ is one-to-one for every $G \in \fD'$; a simple S-scheme is thus always a strong S-scheme. The difference in notation of the homomorphism argument between general S-schemes (in brackets) and simple S-schemes (as subscript) has been chosen in order to facilitate the reading of formulas.

In what follows, we work in parallel with three S-schemes with respect to $\fD'$: a simple S-schemes $\alpha^R$ from $R$ to an auxiliary digraph $\E(R)$, an S-scheme $\rho$ from $R$ to $S$, and a simple S-scheme $\alpha^S$ from $S$ to an auxiliary digraph $\E(S)$. In the case of $\E(R) \in \fD'$, we are in the situation shown in Figure \ref{figure_KommDiagr}. The left triangle $G - \E(S) - S$ and the right triangle $G - \E(R) - R$ commute due to \eqref{defeq_EVsyst}. Because of $\E(R) \in \fD'$ and $\phi_R \in \S( \E(R), R)$, both homomorphisms $\myrhoxi{ \phi_R } : \E(R) \rightarrow S$ and $\myaTxi{S}{\rho( \phi_R )} : \E(R) \rightarrow \E(S)$ are well-defined and strict with $\rho( \phi_R ) = \phi_S \circ \myaTxi{S}{\rho( \phi_R )}$ because of \eqref{defeq_EVsyst}. We are interested in that also the outer triangle $G - \E(R) - S$ commutes
\begin{align} \label{Cond_0}
\forall \; G \in \fD', \xi \in \S(G,R) \mytext{:} \rhoxi & \; = \; \rho( \phi_R ) \circ \aRxi,
\end{align}
which is implied by a commuting inner triangle $G - \E(R) - \E(S)$
\begin{align} \label{Cond_1}
\forall \; G \in \fD', \xi \in \S(G,R) \mytext{:} \myaTxi{S}{\rho(\xi)} & \; = \; \myaTxi{S}{\rho( \phi_R )} \circ \aRxi.
\end{align}

The reason for being interested in \eqref{Cond_0} becomes visible, if we rewrite \eqref{Cond_0}  by means of \eqref{defeq_EVsyst}:
\begin{align*}
\forall \; G \in \fD', \xi \in \S(G,R) \mytext{:} \rhoxi & \; = \; \phi_S \circ \myaTgxi{S}{\E(R)}{\rho( \phi_R )} \circ \aRxi.
\end{align*}
Here, the properties of $\rho$ are mainly determined by the digraphs $\E(R)$ and $\E(S)$ and the strict homomorphism $\eps \equiv \myaTgxi{S}{\E(R)}{\rho( \phi_R )} $ between them. From a theoretical point of view, it is interesting to analyze the mutual dependency between the properties of $\rho$ and $\eps$, and from a practical point of view, there is a perspective to construct an S-scheme $\rho$ from $R$ to $S$ with desired properties by designing (hopefully simple) objects $\E(R)$ and $\E(S)$ and an appropriate homomorphism $\eps$ between them. In order to unburden the notation in what follows, we define for $T \in \{ R, S \}$:
\begin{align*}
\E_o(T) & \; \equiv \; V( \E(T)).
\end{align*}

In the discussion of Figure \ref{figure_KommDiagr}, we have used the assumption $\E(R) \in \fD'$ in order to make the cogwheels interlocking. We give an own name to this assumption and to an additional one:
\begin{itemize}
\item ERD: $\E(R) \in \fD'$;
\item AID: $\myaTgxi{R}{\E(R)}{\phi_R} = \id_{\E_o(R)}$.
\end{itemize}
Assumption AID requires ERD, because otherwise, $\myaTgxi{R}{\E(R)}{\phi_R}$ is not defined. If ERD holds, then also $\rho(\phi_R) : \E(R) \rightarrow S$ is well-defined and strict, because $\rho$ is an S-scheme, and $\myaTgxi{S}{\E(R)}{\rho( \phi_R )} : \E(R) \rightarrow \E(S)$ is strict, too.

In the following definition, we generalize our approach by replacing $\myaTgxi{S}{\E(R)}{\rho( \phi_R )}$ with an arbitrary strict homomorphism from $\E(R)$ to $\E(S)$:

\begin{definition} \label{def_eps}
Let $\eps : \E(R) \rightarrow \E(S)$ be a strict homomorphism. We define for every $\sariegxiStr$
\begin{equation*}
\myegxi{G}{\xi} \quad \equiv \quad \phi_S \circ \eps \circ \aRgxi
\end{equation*}
and call $\eta$ the {\em S-scheme induced by $\eps$.} 

Additionally, we define for all $\sariegxiStr$
\begin{equation*}
E_G(\xi) \quad \equiv \quad \mysetdescr{ v \in G }{ \aSgexiv \in \eps[ \E_o(R) ] }.
\end{equation*}

We write $\exi$ and $E(\xi)$ in the case of a fixed $G \in \fD'$.
\end{definition}

For $\sariegxiStr$, the mapping $\eta_G(\xi) : V(G) \rightarrow V(S)$ is a combination of strict homomorphisms, thus strict. Therefore, $\eta$ is indeed an S-scheme, as suggested by Definition \ref{def_eps}, and $\myaTgxi{S}{G}{\eta(\xi)}$ is well defined. For EDR and $\eps = \myaTgxi{S}{\E(R)}{\rho( \phi_R )}$, we have $\eta(\xi) = \rho( \phi_R ) \circ \aRxi$, hence $\eta = \rho$ in the case of \eqref{Cond_0}. The reader will observe that the set $E_G(\xi)$ can be determined by means of $\eps$ and $\eta_G(\xi)$ for all $\sariegxiStr$; knowledge about $\xi$ is not required.

Induced S-schemes are characterized as follows:
\begin{theorem} \label{theo_induced_exists}
Let $\rho$ be an S-scheme from $R$ to $S$ with respect to $\fD' \subseteq \fD$. If $\rho$ is induced by a strict homomorphism from $\E(R)$ to $\E(S)$ then, for all $G, H \in \fD'$, $\xi \in \S(G,R)$, $\zeta \in \S(H,R)$, $v \in G$, $w \in H$,
\begin{equation} \label{eq_Escheme_rhowert}
\aRgxiv = \myaTgxiv{R}{H}{\zeta}{w} 
\quad \Rightarrow \quad \rhogxiv = \myrhogxiv{H}{\zeta}{w}.
\end{equation}
On the other hand, if ERD and AID are fulfilled and $\rho$ is an S-scheme fulfilling \eqref{eq_Escheme_rhowert}, then $\rho$ is induced by $\myaTgxi{S}{\E(R)}{\rho( \phi_R ) }$.
\end{theorem}
\BP Let $\eps : \E(R) \rightarrow \E(S)$ be a strict homomorphism inducing $\rho$. Then trivially, $\aRgxiv = \myaTgxiv{R}{H}{\zeta}{w}$ yields
\begin{equation*}
\rhogxiv \; = \; \phi_S \left( \eps( \aRgxiv ) \right) \; = \; \phi_S \left( \eps( \myaTgxiv{R}{H}{\zeta}{w} ) \right) \; = \; \myrhogxiv{H}{\zeta}{w},
\end{equation*}
thus \eqref{eq_Escheme_rhowert}.

Now assume ERD and AID and let $\rho$ be an S-scheme fulfilling \eqref{eq_Escheme_rhowert}. Due to AID, $\aRgxiv = \myaTgxiv{R}{\E(R)}{\phi_R}{\aRgxiv}$ for all $\sarierStr$, hence 
\begin{align*}
\rhogxiv
& \; \stackrel{\eqref{eq_Escheme_rhowert}}{=} \; 
\rho_{\E(R)}(\phi_R)(\aRgxiv)
\; \stackrel{\eqref{defeq_EVsyst}}{=} \; 
\phi_S \left( \myaTgxiv{S}{\E(R)}{\rho( \phi_R )}{\aRgxiv} \right),
\end{align*}
and $\rho$ is induced by $\myaTgxi{S}{\E(R)}{\rho( \phi_R ) }$.

\EP

For given $\E(R)$ and $\E(S)$, an induced S-scheme can in general be induced by several homomorphisms $\eps : \E(R) \rightarrow \E(S)$. Nevertheless, the following proposition shows that it is beneficial to select $\eps$ carefully, because suitable properties of $\eps$ guarantee that $\eta$ is ``close to'' being strong:

\begin{proposition} \label{prop_eta_invers}
Let $\eps$ be a strict homomorphism from $\E(R)$ to $\E(S)$. Assume
\begin{equation} \label{epsaepsb_ab}
\forall \; \fa, \fb \in \E_o(R) \mytext{:} \; \; \eps( \fa ) = \eps( \fb ) \; \Rightarrow \; \phi_R( \fa ) = \phi_R( \fb ),
\end{equation}
and assume additionally, that for every $\sariegxiStr$
\begin{equation} \label{condsmall_aexiv_eaxiv}
\forall \; v \in E( \xi ) \; \exists \; \fa \in \myurbild{ \phi_R }( \xiv ) \; \mytext{: } \aSexiv = \eps( \fa ).
\end{equation}
Then, for every $r \in V(R)$ and every $\sariegxiStr$
\begin{equation} \label{eq_eta_invers}
\myurbild{\xi}(r) \cap E( \xi ) \quad = \quad \bigcup_{\fa \in \myurbild{ \phi_R }(r)} \myurbild{ \aSexi }( \eps( \fa ) );
\end{equation}
we can thus reconctruct $\xi \vert_{E(\xi)}$ by means of $\eta(\xi)$ and $\eps$.
\end{proposition}

\BP Let $r \in V(R)$ and let $W$ be the set on the right side of \eqref{eq_eta_invers}. For $v \in \myurbild{\xi}(r) \cap E(\xi)$, assumption \eqref{condsmall_aexiv_eaxiv} delivers an $\fa \in \myurbild{ \phi_R }(r)$ with $\aSexiv = \eps( \fa )$, hence
\begin{equation*}
v \; \in \; \myurbild{ \aSexi }( \aSexiv ) \; = \; \myurbild{ \aSexi }( \eps( \fa ) ).
\end{equation*}
We conclude $v \in W$ due to $\phi_R( \fa ) = r$.

Now let $v \in W$, i.e., $v \in \myurbild{ \aSexi }( \eps( \fa ) )$ for an $\fa \in \myurbild{ \phi_R }(r)$. Then $v \in E(\xi)$. According to \eqref{condsmall_aexiv_eaxiv}, there exists a $\fb \in \myurbild{ \phi_R }( \xiv )$ with $\aSexiv = \eps( \fb )$. Now we get
\begin{equation*}
\eps( \fa ) \; = \; \myaTxiv{S}{\eta(\xi)}{v} \; = \; \eps( \fb ),
\end{equation*}
and assumption \eqref{epsaepsb_ab} yields $r = \phi_R( \fa ) = \phi_R( \fb ) = \xiv$, hence $v \in \myurbild{\xi}(r)$.

\EP

Even if we have found a strict homomorphism $\eps : \E(R) \rightarrow \E(S)$ fulfilling the conditions in Proposition \ref{prop_eta_invers}, there remains a gap for the induced S-scheme $\eta$ to being strong: how to reconstruct $ \xi $ on $ V(G) \setminus E(\xi)$? The gap disappears if $E(\xi) = V(G)$ for all $\sariegxiStr$. The following condition does even more:

\begin{condition} \label{cond_univ_aexiv}
We say that a strict homomorphism $\eps : \E(R) \rightarrow \E(S)$ {\em fulfills Condition \ref{cond_univ_aexiv}}, iff for every $\sariegxiStr$
\begin{equation} \label{cond_aexiv_eaxiv}
\myaTgxi{S}{G}{\eta(\xi)} \quad = \quad \eps \circ \aRgxi.
\end{equation}
\end{condition}

In the case of $\eps = \myaTgxi{S}{\E(R)}{\rho( \phi_R ) }$, this condition is \eqref{Cond_1}. Induced S-schemes fulfilling Condition \ref{cond_univ_aexiv} are characterized in the following theorem:

\begin{theorem} \label{theo_eta_RS} Assume that ERD and AID are fulfilled and that $\rho$ is an S-scheme from $R$ to $S$. We define $\eps \equiv \myaTgxi{S}{\E(R)}{\rho( \phi_R ) }$.

(1) $\rho$ is induced by $ \eps$ and $\eps$ fulfills Condition \ref{cond_univ_aexiv}, iff $\rho$ fulfills the following regularity condition: for every $G, H \in \fD'$, $\xi \in \S(G,R), \zeta \in \S(H,R)$, $v \in V(G)$, $w \in V(H)$
\begin{align} \label{eq_imagebased_gl}
\aRgxiv = \myaTgxiv{R}{H}{\zeta}{w}
& \quad \Rightarrow \quad
\myaTgxiv{S}{G}{\rho(\xi)}{v} = \myaTgxiv{S}{H}{\rho(\zeta)}{w}.
\end{align}
(2) If $\rho$ is induced by $\eps$ and $\eps$ fulfills Condition \ref{cond_univ_aexiv} and \eqref{epsaepsb_ab}, then $\rho$ is a strong S-scheme fulfilling \eqref{eq_imagebased_gl}.

In particular, in the case of $R \in \fD' = \fD$, $R \in \fTa \subseteq \fD' \subseteq \fD$, or $R \in \fP'$ with $\fP' = \fP$ or $\fP' = \fP^*$, $\rho$ can be extended to a strong $\Gamma$-scheme $\rho'$ with $\rho'_G \vert_{\S(G,R)}^{\S(G,S)} = \rho_G$ for all $G \in \fD'$. 
\end{theorem}
\BP (1) If $\rho$ is induced by $ \eps$ and $\eps$ fulfills Condition \ref{cond_univ_aexiv}, then, for $\aRgxiv = \myaTgxiv{R}{H}{\zeta}{w}$,
\begin{equation*}
\myaTgxiv{S}{G}{\rho(\xi)}{v}
\; = \; 
\eps( \aRgxiv )
\; = \; 
\eps( \myaTgxiv{R}{G}{\zeta}{w} )
\; = \;
\myaTgxiv{S}{H}{\rho(\zeta)}{w}.
\end{equation*}

Now assume that the S-scheme $\rho$ fulfills \eqref{eq_imagebased_gl}. According to the second part of Theorem \ref{theo_induced_exists}, $\rho$ is induced by $\eps$, because \eqref{eq_imagebased_gl} implies \eqref{eq_Escheme_rhowert} via \eqref{defeq_EVsyst}. For $\sarierStr$, AID yields $\aRgxiv = \myaTgxiv{R}{\E(R)}{\phi}{\aRgxiv}$, hence,
\begin{equation*}
\myaTgxiv{S}{G}{\rho(\xi)}{v} \; \stackrel{\eqref{eq_imagebased_gl}}{=} \; \myaTgxiv{S}{\E(R)}{\rho(\phi)}{\aRgxiv}) \;  = \; \eps( \aRgxiv ),
\end{equation*}
and $\eps$ fulfills Condition \ref{cond_univ_aexiv}.

(2) Due to part (1) of the theorem, we only have to show that $\rho$ is strong. Because Condition \ref{cond_univ_aexiv} implies $E(\xi) = G$ and \eqref{condsmall_aexiv_eaxiv} for all $\sarierStr$, Proposition \ref{prop_eta_invers} delivers
\begin{equation*}
\myurbild{\xi}(r) \quad = \quad \bigcup_{\fa \in \myurbild{ \phi_R }(r)} \myurbild{ {\myaTgxi{S}{G}{\rho(\xi)}}}( \eps( \fa ) ).
\end{equation*}
for every $\xi \in \S(G,R)$, $G \in \fD'$, $r \in R$. $\rho_G$ is thus one-to-one for every $G \in \fD'$, and $\rho$ is a strong S-scheme. The addendum follows with Proposition \ref{prop_extend_Sscheme}.

\EP

In Theorem \ref{theo_eta_RS_inv} in Section \ref{subsec_inv_Theo_Eta_RS_2}, we will see that with an appropriate choice of $\E(R)$ and $\E(S)$, the inverse of Theorem \ref{theo_eta_RS}(2) is true for the constellations of $R$ and $\fD'$ described in the addendum.

\section{The EV-system of a digraph} \label{sec_EVSystems}

Until now, we have specified nothing about the auxiliary digraphs $\E(R)$ and $\E(S)$ we have used in the definition of an induced S-scheme. In Section \ref{subsec_EVSystems_Props}, we specify $\E(R)$ and $\E(S)$ as {\em EV-systems} of $R$ and $S$, and in Section \ref{subsec_inv_Theo_Eta_RS_2}, we use them in inverting Theorem \ref{theo_eta_RS}(2) for the cases $R \in \fD' = \fD$, $R \in \fTa \subseteq \fD' \subseteq \fD$, and $R \in \fP'$ with $\fP' = \fP$ or $\fP' = \fP^*$. Section \ref{subsec_replacement_Cond1} is devoted to the replacement of Conditon \ref{cond_univ_aexiv} by a more handy condition referring to the homomorphism $\eps : \E(R) \rightarrow \E(S)$ only.

\subsection{EV-systems and their properties} \label{subsec_EVSystems_Props}

In mechanical engineering, the exploded-view drawing of an engine shows the relationship or order of assembly of its components by distributing them in the drawing area in a well-arranged and meaningful way. That is exactly what the {\em EV-system} of a digraph does with respect to the relations between its vertices:

\begin{definition} \label{def_EVsys_alt}
Let $R$ be a digraph and $\fD' \subseteq \fD$. We define
\begin{align*}
\E_o(R) & \equiv \mysetdescr{ ( v, D, U ) }{ v \in V(R), D \subseteq  \Nin_R(v), U \subseteq \Nout_R(v) }.
\end{align*}
For $\fa \in \E_o(R)$, we refer to the three components of $\fa$ by $\fa_1, \fa_2$, and $\fa_3$, and we define
\begin{align*}
\phi_R : \E_o(R) & \rightarrow V(R), \\
\fa & \mapsto \fa_1.
\end{align*}
Furthermore, for every $G \in \fD'$, $\xi \in \S(G,R)$, we define the mapping
\begin{align*}
\myaTgxi{R}{G}{\xi} : V(G) & \rightarrow \E_o(R), \\
v & \mapsto \left( \xiv, \xi[ \Nin_G(v) ], \xi[ \Nout_G(v) ] \right).
\end{align*}
The {\em EV-system $\E(R)$ of $R$ with respect to $\fD'$} is the digraph with vertex set $\E_o(R)$ and arc set $A( \E(R) )$ defined by
\begin{align*}& \fa \fb \in A( \E(R) ) \\
\equiv \quad & \exists \; G \in \fD', \xi \in \S(G,R), vw \in A(G) \mytext{: } \fa = \aRgxiv, \; \fb = \aRgxiw.
\end{align*}
\end{definition}
It is thus the arc set of an EV-system which depends on $\fD'$, whereas the vertex set is independent of it. $\phi_R \circ \aRgxi = \xi$ is trivial. In Lemma \ref{lemma_phi_def}, we will see that $\phi_R$ is strict, as required. As usual, we write $\aRxi$ in the case of a fixed digraph $G$.

\begin{figure}
\begin{center}
\includegraphics[trim = 75 515 270 70, clip]{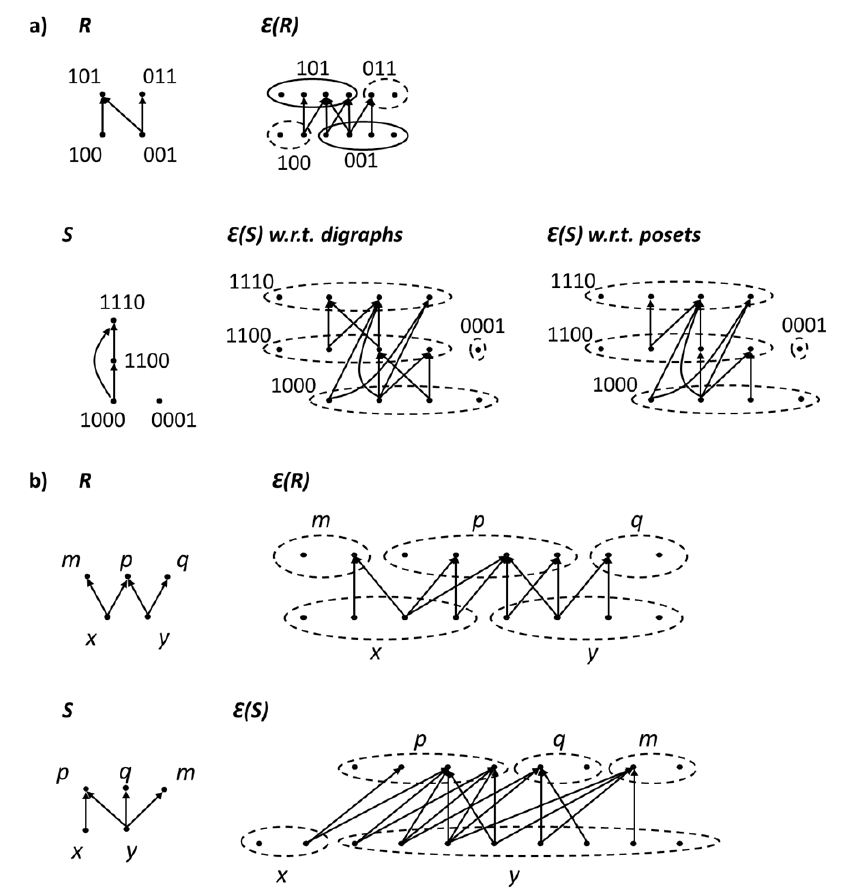}
\caption{\label{figure_Bspl_GKonstr} The posets from Figure \ref{figure_Intro}(a)-(b) drawn as digraphs together with their EV-systems with respect to $\fD$ and $\fP$; all loops are omitted. In the EV-systems, the sets $\myurbild{ \phi_T }(v)$, $v \in V(T)$, $T \in \{ R, S \}$ are encircled and labeled with $v$.}
\end{center}
\end{figure}

Figure \ref{figure_Bspl_GKonstr} shows the EV-systems of the posets in Figure \ref{figure_Intro}(a)-(b) with respect to $\fD$ and $\fP$. All loops are omitted in the diagrams. The poset $S$ in (a) is the only one of the four for which the reference class $\fD'$ makes a difference. The sets $\myurbild{\phi_T}(v)$, $v \in V(T)$, $T \in \{ R, S \}$, are encircled and labeled with the respective $v$. For each point $\fa$ in the diagrams, $\fa_1$ is thus given by this label, and we get $\fa_2$ and $\fa_3$ by looking at the labels of the starting points and end points of arrows ending and starting in $\fa$, respectively.

The reader will have noticed that in the definition of the mappings $\aRgxi$, the restriction to digraphs $G$ {\em contained in $\fD'$} is unnecessary; we could define the mappings in the same way also for {\em all} digraphs $G$. We did not do so for two reasons. Firstly, we would burden the discussion in this section with discriminations about where ERD is required and where not, and secondly, we would gain nothing with this effort, because starting in Section \ref{subsec_inv_Theo_Eta_RS_2}, we work with $\rho$ again and are thus restricted to $G \in \fD'$.

In the following three lemmata, we collect basic properties of the EV-system $\E(R)$ with respect to an arbitrary $\fD' \subseteq \fD$; starting with Definition \ref{def_Xfa}, we deal with the EV-system $\E(R)$ with respect to the classes $\fD'$ we are mainly interested in.

According to the definition of $A( \E(R) )$, the mapping $\agxi$ is a homomorphism from $G$ to $\E(R)$ for every $\sariegxiStr$. For $vw \in A(G^*)$ we have $\xiv \not= \xiw$, hence $\agxiv \not= \agxiw$, and $\agxi$ is strict. $\alpha$ is thus an S-scheme from $R$ to $\E(R)$.

\begin{lemma} \label{lemma_fafb_simpleProp}
For all $\fa, \fb \in \E_o(R)$
\begin{align} \label{fafab_allgemein}
\fa \fb \in A( \E(R ) ) \quad & \Rightarrow \quad \fa_1 \fb_1 \in A(R), \\
\label{fafb_ungleich}
\fa \fb \in A( \E(R)^* ) \quad & \Rightarrow \quad 
\fa_1 \fb_1 \in A( R^* ), \fa_1 \in \fb_2, \fb_1 \in \fa_3, \\
\label{fafb_gleich}
\fa \fb \in A( \E(R) ) \mytext{ with } \fa_1 = \fb_1 \quad & \Rightarrow \quad \fa = \fb. \\
\label{RTa_ERTa}
\mytext{Furthermore,} \quad R \in \fTa \quad & \Rightarrow \quad \E(R) \in \fTa.
\end{align}
\end{lemma}
\BP Let $\fa \fb \in A( \E(R) )$. There exist $G \in \fD', \xi \in \S(G,R)$, and $vw \in A(G)$ with $\fa = \aRxiv$, $\fb = \aRxiw$, thus $\fa_1 \fb_1 = \xiv \xiw \in A(R)$. If $\fa \fb \in A(\E(R)^*)$, then $v \not= w$, hence $\fa_1 \fb_1 = \xi(v) \xi(w) \in A( R^*)$. Furthermore, due to $v \in \Nin_G(w)$, we have
\begin{equation*}
\fa_1 \; = \; \xiv \; \in \; \xi[ \Nin_G(w) ] \; = \; \aRxiw_2 \; = \; \fb_2.
\end{equation*}
$\fb_1 \in \fa_3$ is similarly shown, and \eqref{fafb_ungleich} is proven. For $\fa \fb \in A( \E(R) )$ with $\fa_1 = \fb_1$, \eqref{fafab_allgemein} and \eqref{fafb_ungleich} yield $\fa = \fb$.

Let $R \in \fTa$, and let $\fa^0, \ldots, \fa^I$ be a closed walk in $\E(R)^*$. Due to \eqref{fafb_ungleich}, we have $\fa^{i-1}_1 \fa^i_1 \in A(R^*)$ for all $i \in \myNk{I}$, and the sequence $\fa^0_1, \ldots, \fa^I_1$ is a closed walk in $R^*$, in contradiction to $R \in \fTa$.

\EP

\begin{lemma} \label{lemma_phi_def}
The mapping $\phi_R$ is a strict homomorphism from $\E(R)$ to $R$, and $\alpha^R$ is a simple S-scheme from $R$ to $\E(R)$ with respect to $\fD'$. In the case of ERD,
\begin{align} \label{aphifa_fa}
\myaTxiv{R}{\phi_R}{\fa}_2 \subseteq \fa_2 
& \mytext{ and }
\myaTxiv{R}{\phi_R}{\fa}_3 \subseteq \fa_3
\end{align}
for all $\fa \in \E_o(R)$.
\end{lemma}
\BP $\phi_R$ is a strict homomorphism because of the two first implications in Lemma \ref{lemma_fafb_simpleProp}. We have already seen that $\alpha^R$ is an S-scheme from $R$ to $\E(R)$, and due to $\xi = \phi_R \circ \aRxi$ for all $G \in \fD'$, $\xi \in \S(G,R)$, $\alpha^R$ is simple. 
In the case of ERD, $\alpha_{\phi_R}$ is well defined with
\begin{equation*}
\myaxiv{\phi_R}{\fa}_3 \; = \; \phi_R \left[ \Nout_{\E(R)}( \fa ) \right] \; = \; 
\mysetdescr{ \fb_1 }{ \fb \in \Nout_{\E(R)}( \fa ) }
\; \stackrel{\eqref{fafb_ungleich}}{\subseteq} \; \fa_3.
\end{equation*}
$\myaxiv{\phi_R}{\fa}_2 \subseteq \fa_2$ is shown in the same way.

\EP

\begin{lemma} \label{lemma_phi_fa}
If ERD holds, then AID is equivalent to
\begin{equation} \label{eq_EXI}
\forall\; \fa \in \E_o(R) \; \exists\; G \in \fD', \xi \in \S(G,R), v \in V(G) \mytext{: } \fa = \aRgxiv.
\end{equation}
\end{lemma}
\BP ``$\Rightarrow$'' is trivial. Assume \eqref{eq_EXI}. We write $\phi$ instead of $\phi_R$. Let $\fa \in \E_o(R)$ and $\sarieStr$ with $\fa = \aRgxiv$. Due to \eqref{aphifa_fa}, we have to show $\fa_2 \subseteq \myaTgxiv{R}{\E(R)}{\phi}{\fa}_2$ and $\fa_3 \subseteq \myaTgxiv{R}{\E(R)}{\phi}{\fa}_3$ only.

For $a \in \fa_2$, there exists a $w \in \Nin_G(v)$ with $a = \xi(w)$. $wv \in A(G^*)$ yields $\aRgxiv \in \Nin_{\E(R)}(\fa)$ due to the strictness of $\aRgxi$, hence
\begin{equation*}
a \; = \; \phi( \aRgxiv ) \; \in \; \phi \left[ \Nin_{\E(R)}(\fa) \right] \; = \; \myaTgxiv{R}{\E(R)}{\phi}{\fa}_2.
\end{equation*}
$\fa_3 \subseteq \myaTgxiv{R}{\E(R)}{\phi}{\fa}_3$ is proven in the same way.

\EP

We now show that ERD and AID are fulfilled for the choices of $\fD'$ and $R$ we are particularly interested in: $R \in \fD' = \fD$, $R \in \fTa \subseteq \fD' \subseteq \fD$, $R \in \fD' = \fP$, and $R \in \fD' = \fP^*$. We need the following objects:

\begin{definition} \label{def_Xfa}
For every $m, n \in \myN_0$, we define the digraph $\Xmn \in \fTa$ by
\begin{align*}
V( \Xmn ) & \; \equiv \; D \cup \{ p \} \cup U, \\
A( \Xmn ) & \; \equiv \; \left( D \times \{ p \} \right) \; \cup \; \left( \{ p \} \times U \right).
\end{align*}
where $D$ and $U$ are disjoint sets with $\# D = m$, $\# U = n$, and $p \notin D \cup U$.

Furthermore, for $\fa \in \E_o(R)$, we define the digraph $X(\fa)$
\begin{itemize}
\item as $\myXmnC{\fa_2}{\fa_3}$ in the case of $R \in \fD' = \fD$ or $R \in \fTa \subseteq \fD' \subseteq \fD$;
\item as the transitive hull of $\myXmnC{\fa_2}{\fa_3}$ in the case of $R \in \fD' = \fP^*$;
\item as the transitive hull of $\myXmnC{\fa_2}{\fa_3}$ with loops added for every vertex in the case of $R \in \fD' = \fP$.
\end{itemize}
$\iota(\fa) : V(X(\fa)) \rightarrow V(R)$ is a mapping sending $p$ to $\fa_1$ and $D$ and $U$ bijectively to $\fa_2$ and $\fa_3$, respectively.
\end{definition}

\begin{figure} 
\begin{center}
\includegraphics[trim = 75 680 240 70, clip]{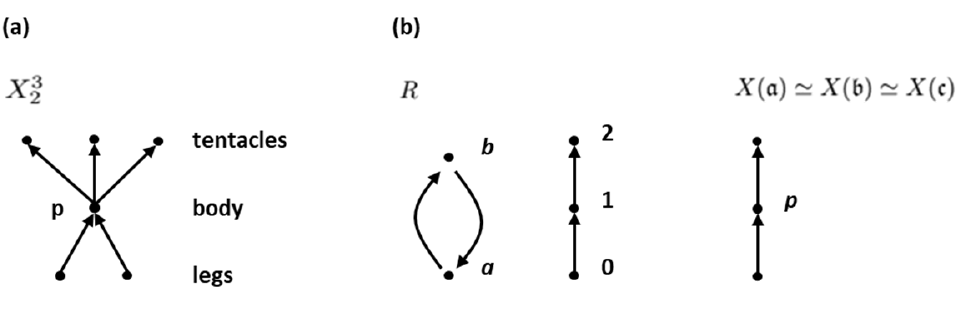}
\caption{\label{figure_Xij} (a) The digraph $\myXmn{2}{3}$ and the nicknames of its vertices.
(b) A digraph $R$ with five vertices. For $\fa \equiv ( a, \{ b \}, \{ b \} )$, $\fb \equiv ( b, \{ a \}, \{ a \} )$, $\fc \equiv ( 1, \{ 0 \},\{ 2 \})$, we have $X(\fa) \simeq X(\fb) \simeq X(\fc)$.}
\end{center}
\end{figure}

$\Xmn$ is thus a bug with $m$ legs, $n$ tentacles and body-vertex $p$. $\myXmn{2}{3}$ is shown in Figure \ref{figure_Xij}(a), and examples for $X(\fa)$ are shown in Figure \ref{figure_Xij}(b). Due to $\fP^*, \fP \subset \fTa$ and $\fP^* \cap \fP = \emptyset$, $X(\fa)$ is in all four cases uniquely determined by $\fD'$, and in all four cases, $X(\fa)$ is an element of $\fD'$. In the following corollary, simple properties of $X(\fa)$ and $\iota(\fa)$ are summarized. It is the inconspicuous first statement which will cause some trouble in Section \ref{sec_undirected} because it does not have a counterpart for undirected graphs.
\begin{corollary} \label{coro_props_Xfa}
Let $R$ and $\fD'$ as in the choices in Definition \ref{def_Xfa}. For every $\fa \in \E_o(R)$, $\pi \in \Aut(X(\fa))$, the vertex $p$ is a fixed point of $\pi$, and $D$ and $U$ are bijectively mapped to $D$ and $U$, respectively. In consequence, for all $\fa \in \E_o(R)$
\begin{align} \label{iotafa_EXI}
\myaxiv{\iota(\fa)}{p} 
& \; = \; \fa, \\ \label{alpha_aiapi}
\myaxiv{\iota(\fa) \circ \pi}{v} 
& \; = \;
\myaiav{\fa}{\pi(v)} \quad \mytext{for all } v \in V(X(\fa)), \pi \in \Aut( X(\fa) ) \\ \label{card_AutXfa}
\# \Aut( X(\fa ) ) & \; = \; ( \# \fa_2 ! ) \cdot ( \# \fa_3 ! ) \; = \; \# \mysetdescr{ \iota(\fa) \circ \pi }{ \pi \in \Aut( X(\fa) ) }.
\end{align}
\end{corollary}
\BP For $\fa_2, \fa_3 \not= \emptyset$, $p$ is the only vertex $v$ of $X(\fa)$ with $\Nin_{X(\fa)}(v)$, $\Nout_{X(\fa)}(v) \not= \emptyset$ (also in the case of $\fa_2 \cap \fa_3 \not= \emptyset$, cf.\ Figure \ref{figure_Xij}(b)). Also in the other cases, $p$ provides a unique empty-nonempty-combination of $\Nin_{X(\fa)}(v)$ and $\Nout_{X(\fa)}(v)$ among the vertices $v$ of $X(\fa)$. Therefore, $p$ is a fixed point of every $\pi \in \Aut( X(\fa) )$, and the rest follows.

\EP

In an intuitive sense, $X(\fa)$ is the most simple digraph $G$ in the respective class $\fD'$ providing a $\xi \in \S(G,R)$ with $\myaxiv{\xi}{p} = \fa $. If $\fa_2$ and $\fa_3$ are not disjoint (which may happen in the case of $ R \in \fD' = \fD$, cf.\ Figure \ref{figure_Xij}(b)), then $\iota(\fa)$ is not one-to-one; however, for our purpose it is enough that it is always strict and fulfills the three equations in Corollary \ref{coro_props_Xfa}.

\begin{figure} 
\begin{center}
\includegraphics[trim = 65 640 320 70, clip]{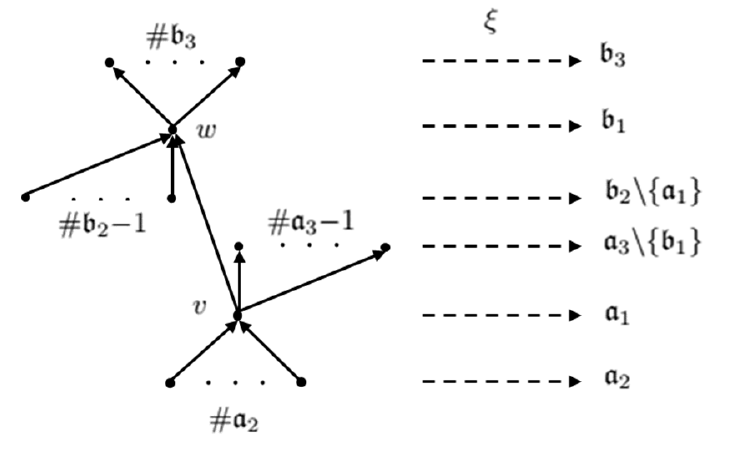}
\caption{\label{figure_H_Gfa_Gfb} The digraph $H \in \fTa$ and the strict homomorphism $\xi : H \rightarrow R$ constructed in the proof of Proposition \ref{prop_fDStr_fTa}.}
\end{center}
\end{figure}

\begin{proposition} \label{prop_fDStr_fTa}
Let $\fTa \subseteq \fD' \subseteq \fD$ and $R \in \fD$. For all $\fa, \fb \in \E_o(R)$,
\begin{align} \label{fafb_ugl_fTa}
\fa \fb \in A( \E(R)^* ) \quad & \Leftrightarrow \quad 
\fa_1 \in \fb_2 \mytext{ and } \fb_1 \in \fa_3, \\ \label{fafb_gl_fTa}
\fa \fa \in A( \E(R) ) \quad & \Leftrightarrow \quad \fa_1 \fa_1 \in A(R).
\end{align}
ERD and AID hold for $\fD' = \fD$ and for $R \in \fTa \subseteq \fD'$.
\end{proposition}
\BP 
Because of Lemma \ref{lemma_fafb_simpleProp}, we have to show ``$\Leftarrow$'' only in \eqref{fafb_ugl_fTa} and \eqref{fafb_gl_fTa}. Let $\fa, \fb \in \E_o(R)$ with $\fa_1 \in \fb_2$ and $\fb_1 \in \fa_3$. Then $\fa_1 \in \Nin_R(\fb_1)$, hence $\fa_1 \fb_1 \in A(R^*)$. We take disjoint isomorphic copies $G_\fa$ and $G_\fb$ of $X(\fa)$ and $X(\fb)$ and connect them to a digraph $H \in \fTa$ as indicated in Figure \ref{figure_H_Gfa_Gfb}: we identify one of the tentacle-vertices of $G_\fa$ with the body-vertex $w$ of $G_\fb$, and we identify one of the leg-vertices of $G_\fb$ with the body-vertex $v$ of $G_\fa$. The mapping $\xi$ from $V(H)$ to $V(R)$ indicated in the figure is a strict homomorphism from $H$ to $R$ with $vw \in A( H^* )$, $\axiv = \fa$,  and $\axiw = \fb$, thus $\fa \fb \in A( \E(R)^* )$.

Let $\fa_1 \fa_1 \in A(R)$. We construct $G \in \fTa$ by adding the loop $(p,p)$ to the arc set of $X(\fa)$, and with \eqref{iotafa_EXI}, we get $\fa \fa = \myagxiv{G}{\iota(\fa)}{p} \, \myagxiv{G}{\iota(\fa)}{p} \in A(\E(R))$.

For $\fD' = \fD$, ERD is trivial, and for $R \in \fTa \subseteq \fD'$, ERD holds due to \eqref{RTa_ERTa}. In both cases, AID follows with \eqref{iotafa_EXI} and Lemma \ref{lemma_phi_fa}.

\EP

\begin{proposition} \label{prop_fDStr_Pos}
Let $\fP'$ be one of the classes $\fP$ or $\fP^*$ and let $R \in \fP'$. Then, for all $\fa, \fb \in \E_o(R)$,
\begin{align} \label{fafb_poset}
\fa \fb \in A( \E(R)^* ) \quad & \Leftrightarrow \quad 
\faaz \in \fb_2 \mytext{ and } \fbbd \in \fa_3,
\end{align}
and ERD and AID hold.
\end{proposition}
\BP Let $\fa \fb \in A( \E(R)^* )$. $\fa_1 \in \fb_2$ and $\fb_1 \in \fa_3$ are due to \eqref{fafb_ungleich}. There exist $P \in \fP'$, $\xi \in \S(P, R)$, and $vw \in A(P^*)$ with $\axiv = \fa$ and $\axiw = \fb$.  The transitivity and antisymmetry of $P$ yield $\Nin_P(v) \subseteq \Nin_P(w)$, hence $\axiv_2 = \xi[ \Nin_P(v) ] \subseteq \xi[ \Nin_P(w) ] = \axiw_2$. The proof of $\fb_3 \subseteq \fa_3$ runs similarly.

For the proof of ``$\Leftarrow$'' in \eqref{fafb_poset}, take the digraph $H$ in Figure \ref{figure_H_Gfa_Gfb} and remove $\# \fb_3$ of the free tentacle-vertices originating from $v$ and $\# \fa_2$ of the free leg-vertices leading to $w$. The transitive hull $P$ of the resulting digraph is an element of $\fP^*$, and for $R \in \fP^*$ and $\xi$ corresponding to the homomorphism in Figure \ref{figure_H_Gfa_Gfb}, we have $\myagxiv{P}{\xi}{v} = \fa$ and $\myagxiv{P}{\xi}{w} = \fb$ with $vw \in A(P^*)$. In the case of $R \in \fP$, add a loop to $A(P)$ for every vertex. \eqref{fafb_poset} is now shown for both choices of $\fP'$.

The transitivity of $\E(R)$ results with \eqref{fafb_poset}. Due to $\fP' \subset \fTa$, we have $\E(R) \in \fTa$ according  to \eqref{RTa_ERTa}, and $\E(R)$ is antisymmetric. For $R \in \fP^*$, the irreflexivity of $\E(R)$ is due to \eqref{fafab_allgemein}. For $R \in \fP$ and $\fa \in \E_o(R)$, we have $\fa \fa = \myagxiv{G}{\iota(\fa)}{p} \, \myagxiv{G}{\iota(\fa)}{p} \in A(\E(R))$ due to \eqref{iotafa_EXI}.

Again, AID follows with \eqref{iotafa_EXI} and Lemma \ref{lemma_phi_fa}.

\EP

\subsection{The inversion of Theorem \ref{theo_eta_RS}(2) } \label{subsec_inv_Theo_Eta_RS_2}

In this section, we prove that Theorem \ref{theo_eta_RS}(2) can be inverted for the choices of $R$ and $\fD'$ we are particularly interested in:

\begin{theorem} \label{theo_eta_RS_inv}
Let $R \in \fD' = \fD$, $R \in \fTa \subseteq \fD' \subseteq \fD$, or $R \in \fD' = \fP'$ with $\fP' = \fP$ or $\fP' = \fP^*$, and let $\E(R)$ and $\E(S)$ be the EV-Systems  of $R$ and $S$. If $\rho$ is a strong S-scheme from $R$ to $S$ fulfilling \eqref{eq_imagebased_gl}, then $\eps \equiv \myaTgxi{S}{\E(R)}{\rho( \phi_R ) }$ fulfills Condition \ref{cond_univ_aexiv}, induces $\rho$, and is one-to-one.
\end{theorem}
The theorem extends also to the sub-classes of $\fD$ mentioned at the end of Section \ref{subsec_DefHomSchemes}.

For the constellations of $R$ and $\fD'$ described in the theorem, ERD and AID are fulfilled according to the Propositions \ref{prop_fDStr_fTa} and \ref{prop_fDStr_Pos}. If $\rho$ fulfills \eqref{eq_imagebased_gl}, then $\eps$ fulfills Condition \ref{cond_univ_aexiv} and induces $\rho$ according to Theorem \ref{theo_eta_RS}(1). What is left to show is that the additional assumption ``$\rho$ strong'' implies ``$\eps$ one-to-one''.

It is beneficial to replace the description $\eps = \myaTgxi{S}{\E(R)}{\rho( \phi_R ) }$ by a more intuitive one. We have 
\begin{equation*}
\myaTxiv{R}{\phi}{\fa}
 \; \stackrel{\mytext{AID}}{=} \; \fa
 \; \stackrel{\eqref{iotafa_EXI}}{=} \; \myaTgxiv{R}{X(\fa)}{\iota(\fa)}{p}
 \end{equation*}
for every $\fa \in \E_o(R)$, hence
\begin{equation} \label{eps_via_aXfaifa}
\eps( \fa) \; = \; 
\myaTgxiv{S}{\E(R)}{\rho( \phi_R ) }{ \fa }
\; \stackrel{\eqref{eq_imagebased_gl}}{=} \; 
\myaTgxiv{S}{X(\fa)}{\rho(\iota(\fa))}{p} \quad \mytext{for all } \fa \in \E(R).
\end{equation}
This is the description of $\eps$ we are using in what follows. The point $p$ and the sets $D$ and $U$ have been specified in the definition of $X(\fa)$.
\begin{corollary} \label{coro_rhophi_DpU}
Let $\fa \in \E_o(R)$. Then, for every $\pi \in \Aut(X(\fa))$,
\begin{align} \label{alpha_ariapi}
\rho(\iota(\fa) \circ \pi)(v) 
& \; = \;
\rho(\iota(\fa))(\pi(v))
\quad \mytext{for all } v \in V(X(\fa)),
\end{align}
and
\begin{align} \label{eqs_riapi}
\begin{split}
\rho(\iota(\fa) \circ \pi)(p) & \; = \; \eps(\fa)_1, \\
\rho(\iota(\fa) \circ \pi)[D] & \; = \; \eps(\fa)_2, \\
\rho(\iota(\fa) \circ \pi)[U] & \; = \; \eps(\fa)_3,
\end{split}
\end{align}
\end{corollary}
\BP Let $\pi \in \Aut(X(\fa))$. According to \eqref{alpha_aiapi}, we have $\myaTxiv{R}{\iota(\fa) \circ \pi}{v} = \myaTxiv{R}{\iota(\fa)}{\pi(v)}$ for all $v \in V(X(\fa))$. Equation \eqref{eq_imagebased_gl} delivers $\myaTxiv{S}{\rho( \iota(\fa) \circ \pi )}{v} = \myaTxiv{S}{\rho(\iota(\fa))}{\pi(v)}$ for all $v \in V(X(\fa))$, and \eqref{alpha_ariapi} follows.

$p$ is a fixed point of $\pi$ according to Corollary \ref{coro_props_Xfa}, hence $\myaTxiv{S}{\rho( \iota(\fa) \circ \pi )}{p} = \myaTxiv{S}{\rho(\iota(\fa))}{p}$, and the first equation follows. Furthermore
\begin{align*}
\rho(\iota(\fa) \circ \pi)[D]
& \; = \;
\mysetdescr{ \rho(\iota(\fa) \circ \pi)(d) }{ d \in D }
\; \stackrel{\eqref{alpha_ariapi}}{=} \;
\mysetdescr{ \rho(\iota(\fa))(\pi(d)) }{ d \in D } \\
& \; = \;
\mysetdescr{ \rho(\iota(\fa))(d) }{ d \in D }
\; = \; 
\rho(\iota(\fa))[D]
\; = \; 
\myaTxiv{S}{\rho(\iota(\fa))}{p}_2
\; \stackrel{\eqref{eps_via_aXfaifa}}{=} \;
\eps( \fa )_2.
\end{align*}
The last equation is proven in the same way.

\EP

\begin{lemma} \label{lemma_rhophi_isom_1}
Let $\fa \in \E_o(R)$. If $\rho$ is strong, then, for every $\pi \in \Aut(X(\fa))$, the mapping $\rho( \iota(\fa) \circ \pi)$ is one-to-one on $D$ and $U$. In particular, $ X( \fa ) \simeq X( \eps( \fa ) )$.
\end{lemma}
\BP Let $\fa \in \E_o(R)$ and $\pi \in \Aut(X(\fa))$. Assume $\rho( \iota(\fa) \circ \pi)(c) = \rho( \iota(\fa) \circ \pi)(d)$ for $c, d \in D$. We define the automorphism $\xi : X(\fa) \rightarrow X(\fa)$ by
\begin{align*}
\xi(v) & \equiv 
\begin{cases}
\pi(v), & \mytext{if} v \in V(X(\fa)) \setminus \{ c, d \}; \\
\pi(d), & \mytext{if} v = c; \\
\pi(c), & \mytext{if} v = d.
\end{cases}
\end{align*}
According to \eqref{alpha_ariapi}, we have $\rho(\iota(\fa) \circ \pi)(v) = \rho(\iota(\fa))(\pi(v)) $
and $\rho(\iota(\fa) \circ \xi)(v) = \rho(\iota(\fa))(\xi(v)) $ for all $v \in V(X(\fa))$. Thus, for every vertex $v \in V( X(\fa) ) \setminus \{ c, d \}$,
\begin{align*}
\myrhoxiv{\iota(\fa) \circ \xi}{v} 
& = \; 
\myrhoxiv{\iota(\fa)}{\xi(v)}
\; = \;
\myrhoxiv{\iota(\fa)}{\pi(v)} \\
& = \;
\myrhoxiv{\iota(\fa) \circ \pi}{v}. \\
\mytext{Furthermore,} \quad \myrhoxiv{\iota(\fa) \circ \xi}{c} & = \; 
\myrhoxiv{\iota(\fa)}{\xi(c)}
\; = \;
\myrhoxiv{\iota(\fa)}{\pi(d)} \\
& = \;
\myrhoxiv{\iota(\fa) \circ \pi}{d}
\; = \;
\myrhoxiv{\iota(\fa) \circ \pi}{c} . \\
\mytext{and similarly} \quad \myrhoxiv{\iota(\fa) \circ \xi}{d} & = \; 
\myrhoxiv{\iota(\fa) \circ \pi}{d},
\end{align*}
hence $\rho(\iota(\fa) \circ \xi) = \rho(\iota(\fa) \circ \pi)$. Because $\rho$ is strong, we have $\iota(\fa) \circ \xi = \iota(\fa) \circ \pi$, thus $c = d$, because all three mappings are one-to-one on $D$. In the same way we see that $\rho( \iota(\fa) \circ \pi)$ is one-to-one on $U$. Now the equations \eqref{eqs_riapi} yield $X(\fa) \simeq X(\eps(\fa))$.

\EP

The following lemma finishes the proof of Theorem \ref{theo_eta_RS_inv}:

\begin{lemma} \label{lemma_eps_oneone}
If $\rho$ is strong, then $\eps$ is one-to-one.
\end{lemma} 
\BP Let $\fa, \fb \in \E(R)$ with $\eps( \fa ) = \eps( \fb )$. According to Lemma \ref{lemma_rhophi_isom_1}, we have $X(\fa) \simeq X(\eps(\fa)) = X(\eps(\fa)) \simeq X(\fb)$. $X(\fa)$ and $X(\fb)$ are thus isomorphic, and due to $X(\fa), X(\fb) \in \fD'$, we have $X(\fa) = X(\fb)$. Let $G \equiv X(\fa)$ and
\begin{align*}
\I(\fa) & \; \equiv \; \mysetdescr{ \iota(\fa) \circ \pi }{\pi \in \Aut(G)}, \\
\I(\fb) & \; \equiv \; \mysetdescr{ \iota(\fb) \circ \pi }{ \pi \in \Aut(G) }.
\end{align*}
With $m \equiv \# \fa_2$, $n \equiv \# \fa_3$, $(m !) \cdot (n !)$ is according to \eqref{card_AutXfa} the cardinality of $\Aut( G )$, $\I(\fa)$, and $\I(\fb)$. $\J \equiv \I(\fa) \cup \I(\fb)$ is a subset of $\S(G, R)$ with $\# \J \geq (m !) \cdot (n !)$; equality holds iff $\fa = \fb$. 

$\rho_G[ \J ]$ is a subset of $\S(G,S)$, and due to Corollary \ref{coro_rhophi_DpU} and Lemma \ref{lemma_rhophi_isom_1}, we have $\# \rho_G[ \J ] = \# \Aut(G) = (m !) \cdot (n !)$ (also in the case of $\eps(\fa)_2 \cap \eps(\fa)_3 \not= \emptyset$). Because $\rho$ is strong, we have $\# \J = \# \rho_G[ \J ]$,  hence $\fa = \fb$.

\EP

Theorem \ref{theo_eta_RS_inv} provides more than the pure inversion of Theorem \ref{theo_eta_RS}(2): it states ``$\eps$ is one-to-one'', whereas in Theorem \ref{theo_eta_RS}(2), the weaker condition ``$\eps(\fa) = \eps( \fb ) \Rightarrow \fa_1 = \fb_1$'' is used. The reason is that in the situation of Theorem \ref{theo_eta_RS_inv}, the extremely simple (and thus: powerful) objects $X(\fa)$ belong to $\fD'$: it are the properties of $X(\fa)$ and $\iota(\fa)$ summarized in Corollary \ref{coro_props_Xfa} which yield the stronger result.

\subsection{The replacement of Condition \ref{cond_univ_aexiv}} \label{subsec_replacement_Cond1}

We have introduced Condition \ref{cond_univ_aexiv} in order to close the gap left by Proposition \ref{prop_eta_invers} to $\eta$ being strong. However, Condition \ref{cond_univ_aexiv} is unwieldy to check because it refers to how $\myagexiv{G}{\xi}{v}$ looks for all $\sarieStr$. It is desirable to have more handy conditions referring to the homomorphism $\eps : \E(R) \rightarrow \E(S)$ only. We need

\begin{lemma} \label{lemma_aexiv_axiv}
Let $R, S \in \fD$, $\fD' \subseteq \fD$. If $\eps : \E(R) \rightarrow \E(S)$ is a strict homomorphism between the EV-systems of $R$ and $S$, then
\begin{align} \label{aexiv_axiv}
\myaTxiv{S}{\eta(\xi)}{v}_2 \subseteq \eps( \aRxiv )_2  & \mytext{ and }
\myaTxiv{S}{\eta(\xi)}{v}_3 \subseteq \eps( \aRxiv )_3,
\end{align}
for all $\sarieStr$, where $\eta(\xi) \in \S(G,S)$ is defined as in Definition \ref{def_eps}. 
\end{lemma}
\BP Let $a \in \aSexiv_2$. There exists a $w \in \Nin_G(v)$ with $\eta(\xi)(w) = a$. We have $\aRxiw \aRxiv \in A( \E(R)^* )$, hence $\eps( \aRxiw ) \eps( \aRxiv ) \in A( \E(S)^* )$, and \eqref{fafb_ungleich} yields $ a = \eta(\xi)(w) = \eps( \aRxiw )_1 \in \eps( \aRxiv )_2$. The second inclusion is shown in the same way.

\EP

Now we can prove

\begin{proposition} \label{prop_repl_Cond1}
Let $R, S \in \fD$, $\fD' \subseteq \fD$, and let $\eps : \E(R) \rightarrow \E(S)$ be a strict homomorphism between the EV-systems of $R$ and $S$. Assume that for $\fa \in \E_o(R)$
\begin{align} \label{cond_epsfa_fa}
\begin{split}
\# \eps( \fa )_2 & \; \leq \; \# \fa_2, \\
\# \eps( \fa )_3 & \; \leq \; \# \fa_3,
\end{split}
\end{align}
and
\begin{align} \label{cond_eps_trennt}
\begin{split}
\forall \; \fb, \fc \in \Nin_{\E(R)}( \fa ) \mytext{: } \eps( \fb )_1 = \eps( \fc )_1 & \; \Rightarrow \; \fb_1 = \fc_1, \\
\forall \; \fb, \fc \in \Nout_{\E(R)}( \fa ) \mytext{: } \eps( \fb )_1 = \eps( \fc )_1 & \; \Rightarrow \; \fb_1 = \fc_1.
\end{split}
\end{align}
Then $\aSexiv = \eps( \aRxiv )$ for all $\sarieStr$ with $\aRxiv = \fa$. In particular, $\eps$ fulfills Condition \ref{cond_univ_aexiv} if \eqref{cond_epsfa_fa} and \eqref{cond_eps_trennt} hold for all $\fa \in \E_o(R)$.
\end{proposition}
\BP Let $\sarieStr$ with $\aRxiv = \fa$. For every $w \in \Nin_G(v)$, we have $\aRxiw \in \Nin_{\E(R)}( \aRxiv )$ due to the definition of $A( \E(R) )$, hence
\begin{align*}
\# \aRxiv_2 & \; \; \; = \;
\# \mysetdescr{ \xi( w ) }{ w \in \Nin_G(v) }
\; \stackrel{\eqref{cond_eps_trennt}}{\leq} \;
\# \mysetdescr{ \eps( \aRxiw )_1 }{ w \in \Nin_G(v) } \\
& \; \; \; = \;
\# \mysetdescr{ \exiw }{ w \in \Nin_G(v) }
\; = \;
\# \aSexiv_2 \\
& \stackrel{\eqref{aexiv_axiv}}{\leq} \;
\# \eps( \aRxiv )_2
\; \stackrel{\eqref{cond_epsfa_fa}}{\leq} \;
\# \aRxiv_2,
\end{align*}
thus $ \# \aSexiv_2 = \# \eps( \aRxiv )_2$. Now \eqref{aexiv_axiv} delivers $ \aSexiv_2 = \eps( \aRxiv )_2 $.
The proof of $ \aSexiv_3 = \eps( \aRxiv )_3 $ is similar, and $\aSexiv = \eps( \aRxiv )$ is shown. The addendum is clear.

\EP

\begin{figure}
\begin{center}
\includegraphics[trim = 70 630 330 70, clip]{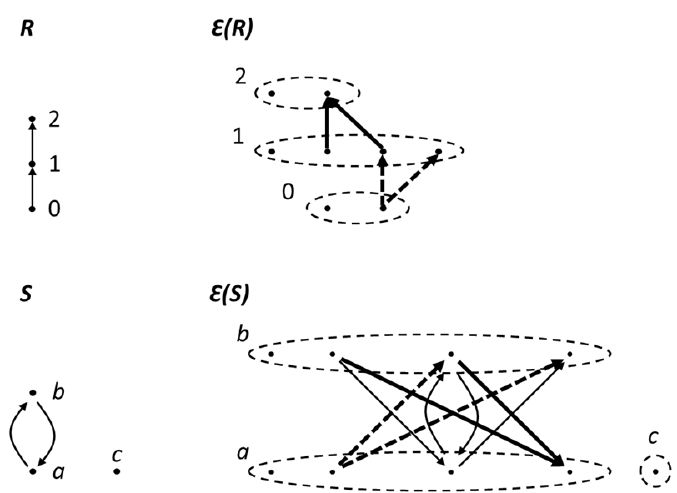}
\caption{\label{figure_RS_nichtRearr} Two digraphs $R$ and $S$ with $R \strG S$ with respect to $\fD$ together with their EV-systems.}
\end{center}
\end{figure}

As an application, Figure \ref{figure_RS_nichtRearr} shows digraphs $R$ and $S$ together with their EV-systems with respect to $\fD$. A homomorphism $\eps$ between the EV-systems is indicated by corresponding bold and dotted shafts of the arrows. The homomorphism fulfills \eqref{cond_epsfa_fa} and \eqref{cond_eps_trennt} for all $\fa \in \E_o(R)$, hence Condition \ref{cond_univ_aexiv}. Because it is additionally one-to-one, it induces a strong S-scheme from $R$ to $S$ with respect to $\fD$  fulfilling \eqref{eq_imagebased_gl} which can be extended to a strong $\Gamma$-scheme from $R$ to $S$ (Theorem \ref{theo_eta_RS}(2)). At the end of Section \ref{subsec_rearr_induced}, we will take the example up again.

Let us examine the assumptions in Proposition \ref{prop_repl_Cond1} more closely! If $\eps$ fulfills Condition \ref{cond_univ_aexiv}, then $\eta$ fulfills \eqref{eq_imagebased_gl} according to Theorem \ref{theo_eta_RS}(1). Looking at the second and third equation in \eqref{eqs_riapi} with $\pi = \id_{V(X(\fa))}$, we conclude that \eqref{cond_epsfa_fa} is necessary for Condition \ref{cond_univ_aexiv} for the choices of $R$ and $\fD'$ we are mainly interested in. However,  \eqref{cond_eps_trennt} is not necessary, as we will show now.

For every flat poset $Q$, we have $\fa_2 = \emptyset$ or $\fa_3 = \emptyset$ for every $\fa \in \E(Q)$, and $\E(Q)$ is a flat poset, too. Let $R$ be a flat connected poset with at least three points, and let $C$ be the two-element chain with $V(C) = \{ 0, 1 \}$ and $A(C) = \{ (0,0), (0,1), (1,1) \}$. Every strict homomorphism $\eps : \E(R) \rightarrow \E(C)$ sends the ``basement'' and the ``upper floor`` of $\E(R)$ to the basement and the upper floor of $\E(C)$, respectively. Following this rule, we define $\eps : \E(R) \rightarrow \E(C)$ by
\begin{align*}
\eps( \fa ) & \equiv 
\begin{cases}
(0, \emptyset, \{ 1 \} ), & \mytext{if } \fa_3 \not= \emptyset \mytext{ (basement to basement)}; \\
(1, \{ 0 \}, \emptyset ), & \mytext{if } \fa_2 \not= \emptyset  \mytext{ (upper floor to upper floor)}; \\
(0, \emptyset, \emptyset), & \mytext{otherwise} \mytext{(isolated points to an isolated point)}.
\end{cases}
\end{align*}
$\eps$ is a strict homomorphism fulfilling \eqref{cond_epsfa_fa} for all $\fa \in \E(R)$. But because $R$ is connected and contains at least three points, $\eps$ violates \eqref{cond_eps_trennt}. However, for all $P \in \fP$, $\xi \in \S(P,R)$, $v \in V(P)$,
\begin{align*}
\aSexiv_2 & = \mysetdescr{ \eps( \aRxiw )_1 }{ w \in \Nin_P(v) } \\
& = 
\begin{cases}
\emptyset, & \mytext{if } \Nin_P(v) = \emptyset; \\
\{ 0 \}, & \mytext{otherwise}
\end{cases} \\
& = \eps( \aRxiv )_2.
\end{align*}
$\aSexiv_3 = \eps( \aRxiv )_3$ is shown in the same way, and $\eps$ fulfills Condition \ref{cond_univ_aexiv}.

\section{The rearrangement method} \label{sec_rearr}

In \cite{aCampo_toappear_0}, the author has developed a method how to rearrange a digraph $R$ in such a way that the relation $R \strG S $ with respect to $\fD$ holds for the digraph $S$ resulting from the rearrangement. In Section \ref{subsec_rearr_induced}, we see that $\rho$ is in fact induced by a strict homomorphism $\eps : \E(R) \rightarrow \E(S)$ between the EV-systems of $R$ and $S$. We describe $\eps$ and analyse its properties. In  Section \ref{subsec_example}, we discuss as examples the pairs of posets in Figure \ref{figure_Intro}(a)-(b) under these view points.

\subsection{{\bf $\rho$} as an induced S-scheme} \label{subsec_rearr_induced}

\begin{figure}
\begin{center}
\includegraphics[trim = 70 640 245 80, clip]{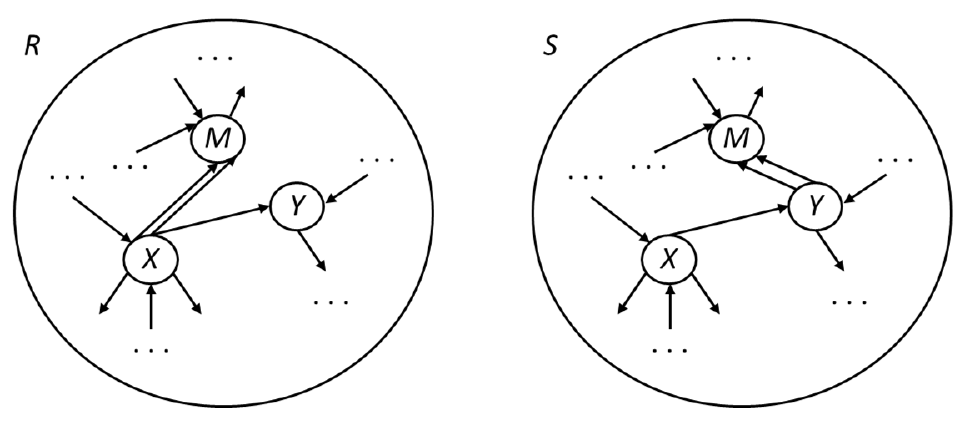}
\caption{\label{fig_abstrConstrMeth} The concept of rearranging a digraph $R$. Explanations in text.}
\end{center}
\end{figure}

The concept of the rearrangement method is illustrated in Figure \ref{fig_abstrConstrMeth}. We have a digraph $R$ and sets $X, Y, M \subseteq V(R)$ with $X \cap M = \emptyset$, $Y \cap M = \emptyset$. We build a new digraph $S$ by replacing all arcs between $M$ and $X$ in $R$ by arcs between $M$ and $Y$. In detail:

\begin{definition}[\cite{aCampo_toappear_0}, Definition 5] \label{def_leq_s}
We agree on the following:
\begin{itemize}
\item $R = (Z, A(R))$ is a digraph.
\item We have disjoint subsets $X$ and $M$ of $Z$.
\item There is a subset $Y \subseteq Z$ with
\begin{align}
\label{bed_WB_empty}
M \cap Y & = \emptyset , \\
\label{bed_WBB_empty}
M \cap N_R(y) & = \emptyset \mytext{ for all } y \in Y,
\end{align}
and $ \beta : X \rightarrow Y $ is a mapping.
\item We define $S$ as the digraph with $V(S) = Z$ and 
\begin{align*}
A(S) & \equiv \; \; A_r \; \cup \; A_d \; \cup \;  A_u \\
\mytext{where} \; A_r & \equiv \; A(R) \setminus \left( ( M \times X ) \cup ( X \times M ) \right), \\
A_d & \equiv \;
\mysetdescr{ m \beta(x) }{ m x \in A(R) \cap ( M \times X )}, \\
A_u & \equiv \;
\mysetdescr{ \beta(x) m }{ x m \in A(R) \cap ( X \times M )}.
\end{align*}
\end{itemize}
\end{definition}

\begin{theorem}[\cite{aCampo_toappear_0}, Theorem 3] \label{theo_rhoxi}
For every $G \in \fD$ and every $\xi \in \S(G,R)$, we define the mapping $\rho_G(\xi) : V(G) \rightarrow Z$ by
\begin{eqnarray} \label{descr_rhoxi}
\forall v \in V(G) \mytext{:} \; \rho_G(\xi)(v) & \equiv & 
\begin{cases}
\beta(\xiv), & \mytext{if } v \in B_\xi, \\
\xiv, & \mytext{otherwise}.
\end{cases}
\end{eqnarray}
where
\begin{equation*}
B_\xi \;  \equiv \; \mysetdescr{ v \in V(G) }{ \xiv \in X \mytext{and} \; \xi[ N_G(v) ] \cap M \not= \emptyset },
\end{equation*}
Assume additionally to the assumptions in Definition \ref{def_leq_s}, that $\beta : R \vert_X \rightarrow R \vert_Y $ is a bijective homomorphism and
\begin{align}
\label{bed_nbh}
\forall x \in X \mytext{:} \; \Nin_R(x) \setminus M \subseteq \Nin_R(  \beta(x) ) \; \mytext{and} \; \Nout_R(x) \setminus M\subseteq \Nout_R( \beta(x) ).
\end{align}
Then $\rho$ is a strong S-scheme from $R$ to $S$ with respect to $\fD$ which can be extended to a strong $\Gamma$-scheme $ \rho'$ with $\rho'_G \vert_{\S(G,R)}^{\S(G,S)} = \rho_G$ for all $G \in \fD$. 
\end{theorem}

It is $\aRgxiv = (\xiv, \xi[ \Nin_G(v) ], \xi[ \Nout_G(v) ])$ which determines $\rhogxiv$ for every $\sarier$. Therefore, \eqref{eq_Escheme_rhowert} holds, and because ERD and AID are fulfilled for $\fD$ and $R$ according to Proposition \ref{prop_fDStr_fTa}, $\rho$ is according to Theorem \ref{theo_induced_exists} induced by
\begin{align*}
\eps : \E(R) & \rightarrow \E(S), \\
\fa & \mapsto \myaTgxiv{S}{\E(R)}{\rho(\phi)}{\fa}.
\end{align*}
We want to describe $\eps( \fa )$ explicitely. We have
\begin{equation*}
B_\phi \quad = \quad \mysetdescr{ \fa \in \E(R)}{ \fa_1 \in X, \phi[ N_{\E(R)}( \fa )] \cap M \not= \emptyset }, 
\end{equation*}
hence, due to \eqref{descr_rhoxi},
\begin{equation} \label{rhophi_rearr}
\eps( \fa )_1 \; = \; \rpha \; =  \; 
\begin{cases}
\beta(\fa_1), & \mytext{if } \fa \in B_\phi, \\
\fa_1, & \mytext{otherwise}.
\end{cases}
\end{equation}

For the determination of $\eps( \fa )_2 = \rho(\phi)\left[ \Nin_{\E(P)}(\fa) \right]$, \eqref{fafb_ugl_fTa} yields
\begin{align*}
\Nin_{\E(P)}(\fa) & \quad = \quad
\mysetdescr{ \fb \in \E(P) }
{\fb_1 \in \fa_2 \; \mytext{and} \; \fa_1 \in \fb_3}.
\end{align*}
Looking at the definition of $B_\phi$ and \eqref{rhophi_rearr}, we see  $\fa_2 \setminus X \subseteq \eps( \fa )_2$. Furthermore,
\begin{itemize}
\item $ \fa_1 \in M$: Then every $\fb \in \Nin_{\E(P)}(\fa) $ with $\fb_1 \in X$ belongs to $B_\phi$ and we conclude $\eps( \fa )_2 = ( \fa_2 \setminus X) \; \cup \; \beta[ \fa_2 \cap X ]$.
\item $ \fa_1 \notin M$: Let $x \in \fa_2 \cap X$. We have $\left( x, \emptyset, \{ \fa_1 \} \right) \in \left( \Nin_{\E(P)}(\fa) \right) \setminus B_\phi$, thus $x \in \eps( \fa )_2$. Additionally, under all elements of $\Nin_{\E(P)}(\fa)$ with first component $x$, it is $(x, \Nin_R(x), \Nout_R(x))$ which has the largest second and third component; therefore, 
$\beta( x ) \in \eps( \fa )_2$, if the intersection $N_R(x) \cap M$ is not empty.
\end{itemize}
Making the same considerations for $\eps( \fa )_3$, we get all together
\begin{align} \label{eps_rearr}
\begin{split}
\eps( \fa )_1 & = \rpha, \\
\eps( \fa )_2 & =
\begin{cases}
( \fa_2 \setminus X ) \; \cup \; \beta[ \fa_2 \cap X ] \quad \quad \quad \quad \quad \quad \quad \quad \quad \;  
\mytext{if } \fa_1 \in  M; & \\
\fa_2 \; \cup \; \beta \left[ \mysetdescr{ x \in \fa_2 \cap X }{ N_R(x) \cap M \not= \emptyset } \right] \quad \mytext{otherwise}; &
\end{cases} \\
\eps( \fa )_3 & =
\begin{cases}
( \fa_3 \setminus X ) \; \cup \; \beta[ \fa_3 \cap X ] \quad \quad \quad \quad \quad \quad \quad \quad \quad \; 
\mytext{if } \fa_1 \in  M; & \\
\fa_3 \; \cup \; \beta \left[ \mysetdescr{ x \in \fa_3 \cap X }{ N_R(x) \cap M \not= \emptyset  \right]} \quad \mytext{otherwise}. &
\end{cases}
\end{split}
\end{align}

Because $\rho$ is induced and strong, Theorem \ref{theo_eta_RS}(1) and Theorem \ref{theo_eta_RS_inv} yield
\begin{align} \label{eps_Cond1_1to1}
\begin{split}
& \eps \; \mytext{fulfills Condition \ref{cond_univ_aexiv}} \\
\Leftrightarrow \quad &
\rho \; \mytext{fulfills } \eqref{eq_imagebased_gl} \\
\Rightarrow \quad & \eps \; \mytext{is one-to-one.} 
\end{split}
\end{align}
But in many cases, the homomorphism $\eps$ resulting from the rearrangement method will not be one-to-one, because, in the case of $X \cap Y = \emptyset$, the following lemma states that $\eps$ will be one-to-one iff every non-isolated $x \in X$ is either encapsulated by $M$ or by $V(R) \setminus M$:

\begin{lemma} \label{lemma_eps_nicht_1t1_Cond1}
Let $R$, $S$, and $\rho$ as in Theorem \ref{theo_rhoxi} and $\eps \equiv \alpha_{\rho(\phi)}$. If 
\begin{align} \label{bed_eps_1t1}
\forall \; x \in X \mytext{: } N_R(x) \cap M = \emptyset & \mytext{ or } N_R(x) \setminus M = \emptyset.
\end{align}
then $\eps$ is one-to-one. In the case of $X \cap Y = \emptyset$, the inverse is true, too.
\end{lemma}
\BP ``$\Rightarrow$'': If \eqref{bed_eps_1t1} holds, then the equations \eqref{eps_rearr} become
\begin{align*}
\eps( \fa )_1 & = \rpha, \\
\eps( \fa )_2 & =
\begin{cases}
( \fa_2 \setminus X ) \; \cup \; \beta[ \fa_2 \cap X ] & \mytext{if } \fa_1 \in  M; \\
\; \fa_2 & \mytext{otherwise};
\end{cases} \\
\eps( \fa )_3 & =
\begin{cases}
( \fa_3 \setminus X ) \; \cup \; \beta[ \fa_3 \cap X ] & \mytext{if } \fa_1 \in  M; \\
\; \fa_3 & \mytext{otherwise}.
\end{cases}
\end{align*}
Let $\fa \in \E_o(R)$. If $\eps( \fa )_1 \in M$, then $\fa_1 = \eps( \fa )_1$ due to \eqref{bed_WB_empty}, and $\fa_2 \cap Y = \emptyset = \fa_3 \cap Y$ due to \eqref{bed_nbh}. Therefore,
\begin{align*}
\fa_2 & \; = \;
\left( \eps( \fa )_2 \setminus Y \right) \; \cup \; \myurbild{\beta}\left( \eps( \fa )_2 \cap Y \right), \\
\fa_3 & \; = \;
\left( \eps( \fa )_3 \setminus Y \right) \; \cup \; \myurbild{\beta}\left( \eps( \fa )_3 \cap Y \right).
\end{align*}
If $\eps( \fa )_1 \notin M$, then $\fa_1 \notin M$ and $\eps(\fa)_2 = \fa_2$, $\eps(\fa)_3 = \fa_3$. If the set $( \fa_2 \cup \fa_3 ) \cap M$ is empty, then the triplet $\fa$ cannot be an element of $B_\phi$, hence $\fa_1 = \eps( \fa_1 )$. And in the case of $( \fa_2 \cup \fa_3 ) \cap M \not= \emptyset$, we have $\fa_1 = \myurbild{\beta}( \eps( \fa )_1 )$ for $\eps( \fa )_1 \in Y$ (use \eqref{bed_nbh}) and $\fa_1 = \eps( \fa )_1$ for $\eps( \fa )_1 \notin Y$.

``$\Leftarrow$'': Assume that \eqref{bed_eps_1t1} does not hold for $x \in X$. Select a vertex $v \in N_R(x) \setminus M$. In the case of $v \in \Nin_R(x)$, define
\begin{align*}
\fa & \; \equiv \; ( v, \emptyset, \{ x \} ), \\
\fb & \; \equiv \; ( v, \emptyset, \{ x, \beta(x) \} ).
\end{align*}
Then $\fa \in \E_o(R)$, and due to $v \in \Nin_R( x ) \setminus M \stackrel{\eqref{bed_nbh}}{\subseteq} \Nin_R( \beta(x) )$, we also have $\fb \in \E_o(R)$. Due to $X \cap Y = \emptyset$, we have $\fa \not= \fb$, but the formulas in \eqref{eps_rearr} yield $\eps( \fa ) = ( v, \emptyset, \{ x, \beta(x)\} ) = \eps( \fb )$ (for the latter equality, we need $X \cap Y = \emptyset$ again). In the case of $v \in \Nout_R(x)$, work with $\fa \equiv ( v, \{ x \}, \emptyset )$, $ \fb \equiv ( v, \{ x, \beta(x)\}, \emptyset )$.

\EP

The condition $X \cap Y = \emptyset$ cannot easily be skipped in this Lemma. For $X = Y$ and $\beta$ being the identity mapping of $X$, the rearrangement method delivers $S = R$, $\rho$ is the trivial $\Gamma$-scheme with $\rho_G = \id_{\H(G,R)}$ for all $G \in \fD$, and $\eps$ is the identity mapping of $\E_o(R)$ according to \eqref{eps_rearr} and \eqref{rhophi_rearr}. $\eps$ is thus one-to-one, whatever the structure of $R$ is.

We want to show that the digraph $S$ in Figure \ref{figure_RS_nichtRearr} in Section \ref{subsec_replacement_Cond1} cannot be constructed by the rearrangement of the digraph $R$ in the figure. Due to \eqref{bed_WB_empty} and \eqref{bed_WBB_empty}, neither the set $M$ nor the set $Y$ can contain the vertex $1$. These sets must be singletons, one of them containing the vertex $0$, the other one the vertex $2$. If the set $X$ contains $1$, the rearrangement of $R$ results in a digraph with V-shaped or $\Lambda$-shaped diagram (and Theorem \ref{theo_rhoxi} cannot be applied because \eqref{bed_nbh} is violated). And in the case of $1 \notin X$, we have $X = \emptyset$ or $X = Y$, and the rearrangement method produces nothing than $R$.

On the one hand, we have thus a rearrangement method producing induced strong S-schemes $\rho$ for which $\eps$ does not fulfill Condition \ref{cond_univ_aexiv} in many cases, and on the other hand, we have Theorem \ref{theo_eta_RS}(2) which states that many homomorphisms $\eps$ fulfilling Condition \ref{cond_univ_aexiv} induce a strong S-scheme $\rho$. There is still something to do in characterizing digraphs $R$ and $S$ with a strong induced S-scheme between them.

\subsection{Examples} \label{subsec_example}

\begin{figure}
\begin{center}
\includegraphics[trim = 70 550 210 60, clip]{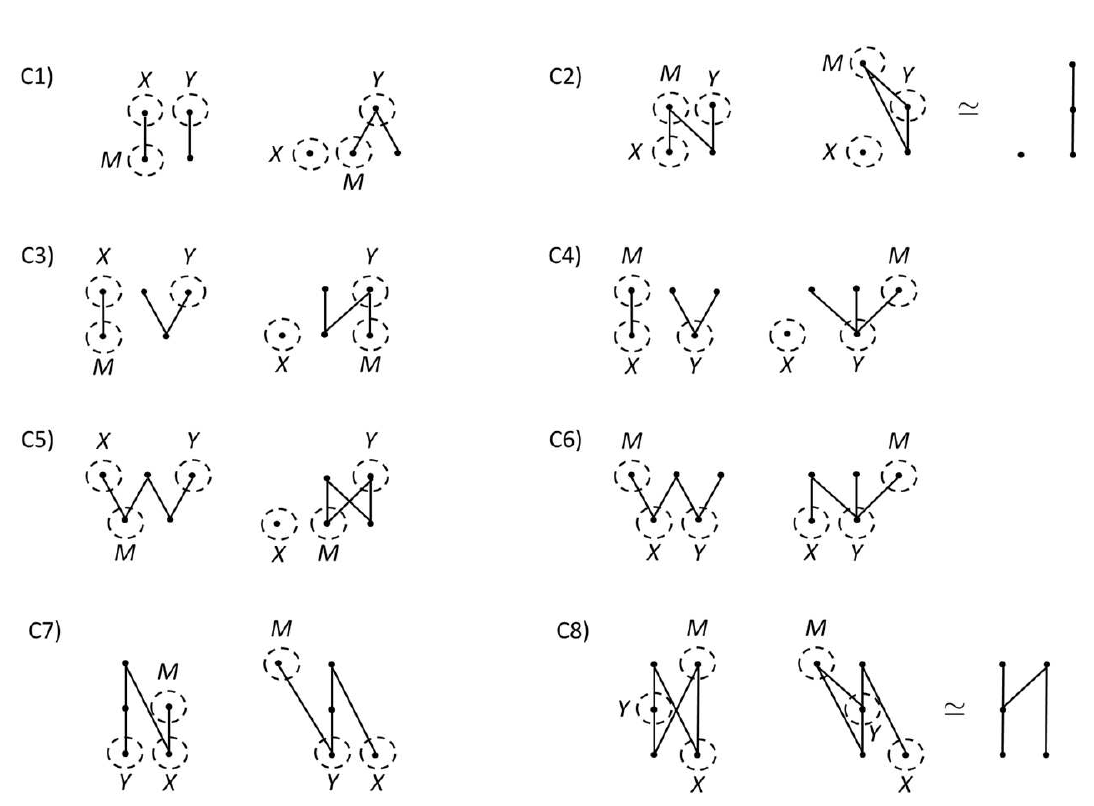}
\caption{\label{fig_TableConstr} Hasse-diagrams of eight pairs of posets $R$ and $S$ for which $R \strG S$ has been shown in \cite{aCampo_toappear_0} by means of the rearrangement method. In all examples, $R$ is on the left and $S$ on the right.}
\end{center}
\end{figure}

Figure \ref{fig_TableConstr} shows the Hasse-diagrams of eight pairs of posets $R$ and $S$ for which $R \strG S$ has been shown in \cite{aCampo_toappear_0} by means of the rearrangement method described in the previous section. The respective one-element vertex sets $X$, $Y$, and $M$ are marked. In all cases, the respective strong S-scheme $\rho$ from $R$ to $S$ is induced by the strict homomorphism $\eps : \E(R) \rightarrow \E(S)$ described by the equations \eqref{eps_rearr}.

\begin{figure}[h]
\begin{center}
\includegraphics[trim = 70 495 190 70, clip]{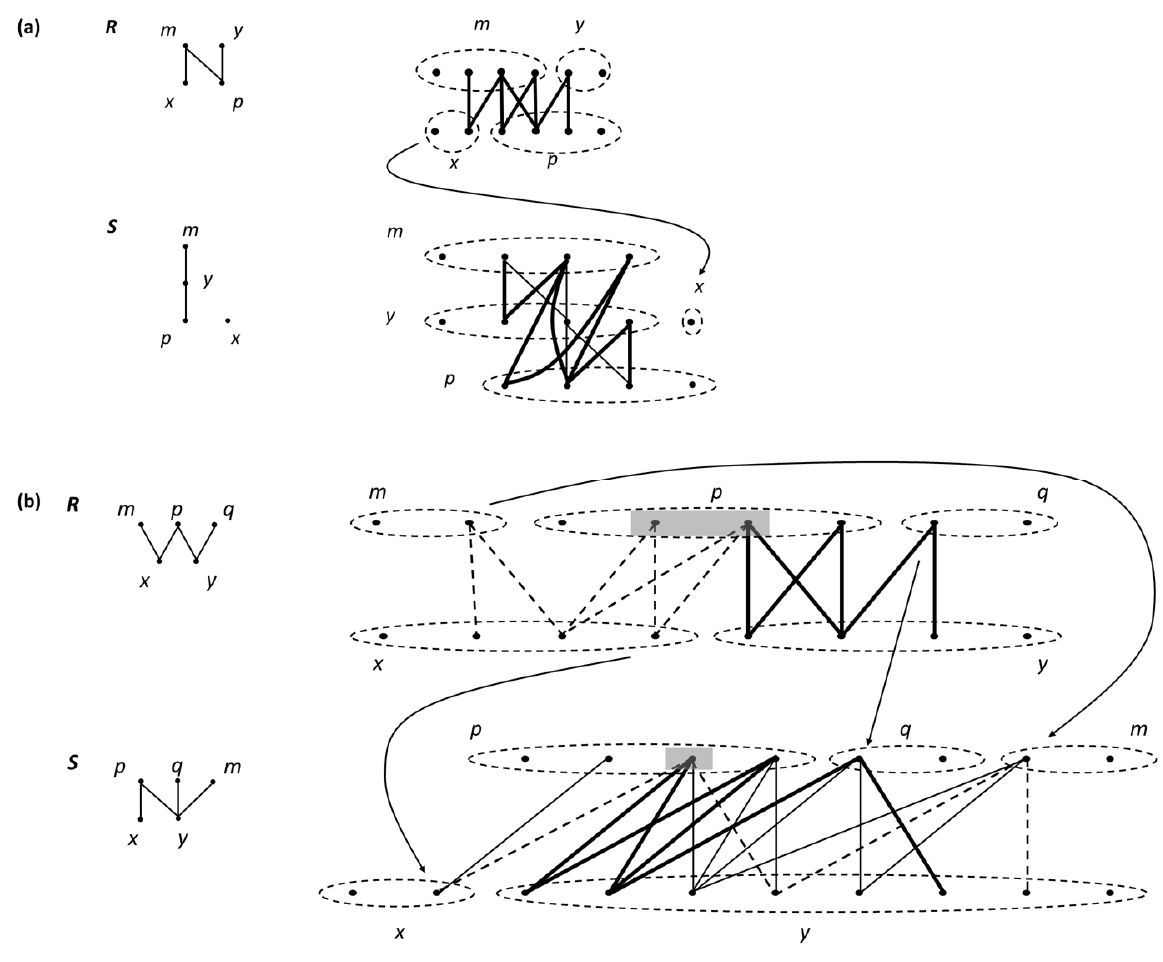}
\caption{\label{figure_Bspl_C7} Hasse-diagrams of the posets $R$ and $S$ from Figure \ref{figure_Intro}(a)-(b), their EV-systems, and a homomorphisms between them. Explanations in text.}
\end{center}
\end{figure}

Figure \ref{figure_Bspl_C7} shows the Hasse-diagrams of the posets $R$ and $S$ from Figure \ref{figure_Intro}(a)-(b), their EV-systems, and the respective strict homomorphism $\eps$ between them; the pairs are the second and the sixth one in Figure \ref{fig_TableConstr}.

In Figure \ref{figure_Bspl_C7}(a), the main part of $\eps$ is indicated by the bold lines in the EV-systems. $\eps$ is one-to-one and by means of Proposition \ref{prop_repl_Cond1}, it is easily seen that it fulfills Condition \ref{cond_univ_aexiv}, too. According to Theorem \ref{theo_eta_RS}, $\rho$ fulfills \eqref{eq_imagebased_gl}. We have 
\begin{equation*}
B_\phi \; = \; \{ ( x, \emptyset, \{ m \} )\}
\end{equation*}
and the total homomorphism is given by
\newline

\parbox{5cm}{
\begin{tabular}{| l | l | }
\hline
$ \fa \in \E_o(R) $ & $\eps( \fa ) \in \E_o(S)$ \\
\hline \hline
$ ( x, \emptyset, \emptyset ) $ &
$ ( x, \emptyset, \emptyset ) $ \\
$ ( x, \emptyset, \{ m \} ) $ &
$ ( y, \emptyset, \{ m \} ) $ \\
\hline
$ ( p, \emptyset, \{ m \} ) $ &
$ ( p, \emptyset, \{ m \} ) $ \\
$ ( p, \emptyset, \{ m, y \} ) $ &
$ ( p, \emptyset, \{ m, y \} ) $ \\
$ ( p, \emptyset, \{ m \} ) $ &
$ ( p, \emptyset, \{ m \} ) $ \\
$ ( p, \emptyset, \emptyset ) $ &
$ ( p, \emptyset, \emptyset ) $ \\
\hline
\end{tabular} }
\parbox{5cm}{
\begin{tabular}{| l | l | }
\hline
$ \fa \in \E_o(R) $ & $\eps( \fa ) \in \E_o(S)$ \\
\hline \hline
$ ( m, \emptyset, \emptyset ) $ & 
$ ( m, \emptyset, \emptyset ) $ \\
$ ( m, \{ x \}, \emptyset ) $ & 
$ ( m, \{ y \}, \emptyset ) $ \\
$ ( m, \{ x, p \}, \emptyset ) $ & 
$ ( m, \{ y, p \}, \emptyset ) $ \\
$ ( m, \{ p \}, \emptyset ) $ & 
$ ( m, \{ p \}, \emptyset ) $ \\
\hline
$ ( y, \{ p \}, \emptyset ) $ & 
$ ( y, \{ p \}, \emptyset ) $ \\
$ ( y, \emptyset, \emptyset ) $ & 
$ ( y, \emptyset, \emptyset ) $ \\
\hline
\end{tabular} }
\newline
\newline

The pairs of posets (C1) and (C3)-(C5) in Figure \ref{fig_TableConstr} have the same properties: the respective strong S-scheme $\rho$ fulfills \eqref{eq_imagebased_gl} and is induced by a one-to-one homomorphism from $\E(R)$ to $\E(S)$ fulfilling Condition \ref{cond_univ_aexiv}.

For the example in Figure \ref{figure_Bspl_C7}(b), three point mappings are indicated by arrows. The rest of $\eps$ is easily seen as follows: The sub-poset of $R$ drawn with bold lines is mapped to the isomorphic bold-lined sub-poset of $S$, whereas the sub-poset of $R$ drawn with dotted lines is flipped and mapped to the M-shaped dotted sub-poset of $S$; the two shaded points $( p, \{ x \}, \emptyset )$ and $( p, \{ x, y \}, \emptyset )$ of $\E(R)$ are both mapped to the shaded point $( p, \{ x, y \}, \emptyset )$ of $\E(S)$. According to \eqref{eps_Cond1_1to1}, $\eps$ cannot fulfill Condition \ref{cond_univ_aexiv} and $\rho$ cannot fulfill \eqref{eq_imagebased_gl}. The total homomorphism is listed in the following table; we have 
\begin{equation*}
B_\phi \; = \; \{ ( x, \emptyset, \{ m \} ), ( x, \emptyset, \{ m, p \} ) \}.
\end{equation*}
\newline
\parbox{5cm}{
\begin{tabular}{| l | l | }
\hline
$ \fa \in \E_o(R) $ & $\eps( \fa ) \in \E_o(S)$ \\
\hline \hline
$ ( x, \emptyset, \emptyset ) $ &
$ ( x, \emptyset, \emptyset ) $ \\
$ ( x, \emptyset, \{ m \} ) $ &
$ ( y, \emptyset, \{ m \} ) $ \\
$ ( x, \emptyset, \{ m, p \} ) $ &
$ ( y, \emptyset, \{ p, m \} ) $ \\
$ ( x, \emptyset, \{ p \} ) $ &
$ ( x, \emptyset, \{ p \} ) $ \\
\hline
$ ( y, \emptyset, \emptyset ) $ &
$ ( y, \emptyset, \emptyset ) $ \\
$ ( y, \emptyset, \{ q \} ) $ &
$ ( y, \emptyset, \{ q \} ) $ \\
$ ( y, \emptyset, \{ p, q \} ) $ &
$ ( y, \emptyset, \{ p, q \} ) $ \\
$ ( y, \emptyset, \{ p \} ) $ &
$ ( y, \emptyset, \{ p \} ) $ \\
\hline
\end{tabular} }
\parbox{5cm}{
\begin{tabular}{| l | l | }
\hline
$ \fa \in \E_o(R) $ & $\eps( \fa ) \in \E_o(S)$ \\
\hline \hline
$ ( m, \emptyset, \emptyset ) $ & 
$ ( m, \emptyset, \emptyset ) $ \\
$ ( m, \{ x \}, \emptyset ) $ & 
$ ( m, \{ y \}, \emptyset ) $ \\
\hline
$ ( p, \emptyset, \emptyset ) $ & 
$ ( p, \emptyset, \emptyset ) $ \\
$ ( p, \{ x \}, \emptyset ) $ & 
$ ( p, \{ x, y \}, \emptyset ) $ \\
$ ( p, \{ x, y \}, \emptyset ) $ & 
$ ( p, \{ x, y \}, \emptyset ) $ \\
$ ( p, \{ y \}, \emptyset ) $ & 
$ ( p, \{ y \}, \emptyset ) $ \\
\hline
$ ( q, \emptyset, \emptyset ) $ & 
$ ( q, \emptyset, \emptyset ) $ \\
$ ( q, \{ y \}, \emptyset ) $ & 
$ ( q, \{ y \}, \emptyset ) $ \\
\hline
\end{tabular} }
\newline
\newline

Let $\fc \equiv ( p, \{x\}, \emptyset ) \in \E_o(R)$ and $\fd \equiv ( p, \{ x, y \}, \emptyset ) \in \E_o(R)$. For all points $\fa \in \E_o(R) \setminus \{ \fc, \fd \}$, the conditions in Proposition \ref{prop_repl_Cond1} are fulfilled, but $\fc$ violates \eqref{cond_epsfa_fa} and $\fd$ violates \eqref{cond_eps_trennt}. The points $\fc$ and $\fd$ are thus the only points $\fa \in \E_o(R)$ for which $\sarie$ may exist with $\aRxiv = \fa $ and $\aSexiv \not= \eps( \fa )$.

\begin{figure}
\begin{center}
\includegraphics[trim = 70 520 280 115, clip]{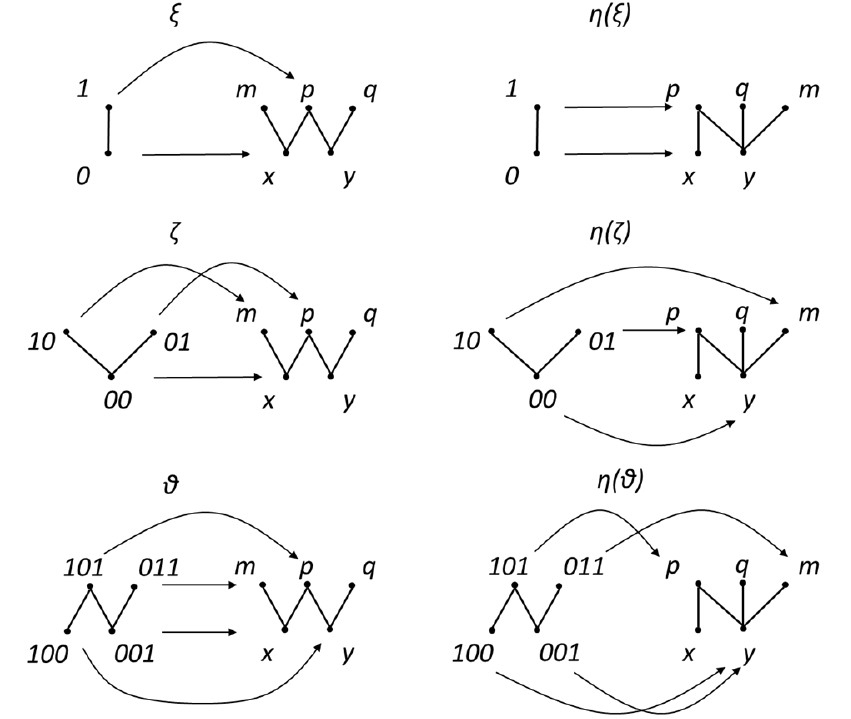}
\caption{\label{figure_BsplHomsC7} Homomorphisms $\xi \in \S( C, R)$, $\zeta \in \S( V, R)$, and $\theta \in \S( N, R)$, and their images under $\eta$.}
\end{center}
\end{figure}

Figure \ref{figure_BsplHomsC7} shows the Hasse-diagrams of three posets $C$, $V$, and $N$, and for each of them a strict homomorphism to $R$ and its image resulting under $\eta$. $C$ is the two-element chain defined in Section \ref{subsec_replacement_Cond1}, and $V$ and $N$ are the posets with V- and N-shaped Hasse-diagrams. In the following tables, the objects of interest in $\E(R)$ and $\E(S)$ are listed.
\newline

\begin{tabular}{| c | l | l | l | }
\hline
$ v \in V(C) $ & $ \axiv $ & $\eps( \axiv ))$ & $ \myaexiv{\xi}{v} $\\
\hline \hline
0 & $( x, \emptyset, \{ p \} )$ & $( x, \emptyset, \{ p \} )$ & $( x, \emptyset, \{ p \} )$ \\
1 & $( p, \{ x \}, \emptyset )$ & $( p, \{ x, y \}, \emptyset )$ & $( p, \{ x \}, \emptyset )$ \\
\hline
\end{tabular}
\newline
\newline

\begin{tabular}{| c | l | l | l | }
\hline
$ v \in V(V) $ & $ \myaxiv{\zeta}{x} $ & $\eps( \myaxiv{\zeta}{x} )$ & $ \myaexiv{\zeta}{x} $ \\
\hline \hline
00 & $( x, \emptyset, \{ m, p \} )$ & $( y, \emptyset, \{ p, m\} )$ & $( y, \emptyset, \{ p, m\} )$ \\
10 & $( m, \{ x \}, \emptyset )$ & $( m, \{ y \}, \emptyset )$ & $( m, \{ y \}, \emptyset )$ \\
01 & $( p, \{ x \}, \emptyset )$ & $( p, \{ x, y \}, \emptyset )$ & $( p, \{ y \}, \emptyset )$ \\
\hline
\end{tabular}
\newline
\newline

\begin{tabular}{| c | l | l | l | }
\hline
$ v \in V(N) $ & $ \myaxiv{\theta}{x} $ & $\eps( \myaxiv{\theta}{x} )$ & $ \myaexiv{\theta}{x} $ \\
\hline \hline
100 & $( y, \emptyset, \{ p \} )$ & $( y, \emptyset, \{ p \} )$ & $( y, \emptyset, \{ p \} )$ \\
001 & $( x, \emptyset, \{ m, p \} )$ & $( y, \emptyset, \{ p, m \} )$ & $( y, \emptyset, \{ p, m \} )$ \\
101 & $( p, \{ x, y \}, \emptyset )$ & $( p, \{ x, y \}, \emptyset )$ & $( p, \{ y \}, \emptyset )$ \\
011 & $( m, \{ x \}, \emptyset )$ & $( m, \{ y \}, \emptyset )$ & $( m, \{ y \}, \emptyset )$ \\
\hline
\end{tabular}
\newline
\newline

We have thus $ \myaxiv{\xi}{1} = \myaxiv{\zeta}{01} = \fc$, but 
$ \eps( \fc ) \not= \myaexiv{\xi}{1} \not= \myaexiv{\zeta}{01} \not= \eps( \fc )$, and we have $ \myaxiv{\theta}{101} = \fd $, but $ \myaexiv{\theta}{101} \not= \eps( \fd )$. 

At the end of this section, we show that no strong S-scheme $\rho$ exists from $R$ to $S$ fulfilling \eqref{eq_imagebased_gl}. Assume that such a strong S-scheme exists. According to Theorem \ref{theo_eta_RS_inv}, $\rho$ is induced by $\eps \equiv \alpha_{\rho(\phi)}$. Now let $P \equiv R$ and let $\xi \equiv \id_{V(R)}$ be the identity mapping of $R$. For $\myaxiv{\xi}{x} = ( x, \emptyset, \{ m, p \} )$ and $\myaxiv{\xi}{y} = ( y, \emptyset, \{ p, q \} )$, Lemma \ref{lemma_rhophi_isom_1} delivers $\# \eps( \myaxiv{\xi}{x} )_3 = 2 = \# \eps( \myaxiv{\xi}{y} )_3$. Because $y$ is the only point $v \in V(S)$ with $\# \Nout_S(v) > 1$, we conclude
\begin{equation*}
\eta(\xi)(x) \; = \; \eps( \myaxiv{\xi}{x} )_1 \; = \; y
 \; = \; \eps( \myaxiv{\xi}{y} )_1 \; = \; \eta(\xi)(y).
\end{equation*}
Therefore,
\begin{equation*}
\myaexiv{\xi}{p}_2 \; = \; \eta(\xi)[ \Nin_G(p) ] \; = \; \eta(\xi)[ \{ x, y \} ] \; = \; \{ y \}.
\end{equation*}
But for $\myaxiv{\xi}{p} = (p, \{x,y\},\emptyset)$, Lemma \ref{lemma_rhophi_isom_1} delivers $\# \eps( \myaxiv{\xi}{p} )_2 = 2$. We have thus $ \myaexiv{\xi}{p}_2 \not= \eps( \myaxiv{\xi}{p} )_2$, and $\eps$ does not fulfill Condition \ref{cond_univ_aexiv}. But due to Theorem \ref{theo_eta_RS_inv}, this is a contradiction to $\rho$ being a strong S-scheme fulfilling \eqref{eq_imagebased_gl}.

Also for the pairs (C7) and (C8) in Figure \ref{fig_TableConstr}, the respective $\eps$ is not one-to-one and does not fulfill Condition \ref{cond_univ_aexiv}, as is easily seen by means of Lemma \ref{lemma_eps_nicht_1t1_Cond1}. In consequence, the respective S-scheme does not fulfill \eqref{eq_imagebased_gl}. For the pair $R$, $S$ in (C7), we see by inspection of $\id_{V(R)}$ and the consequences for $\eta( \id_{V(R)} )$ that there exists no strong S-scheme from $R$ to $S$ fulfilling \eqref{eq_imagebased_gl}.

\section{Undirected graphs} \label{sec_undirected}

Let $\fS$ be a representative system of the non-empty finite {\em symmetric} digraphs and $\fU$ a representative system of the non-empty finite {\em undirected} graphs. Already in \cite[Section 4.2]{aCampo_toappear_0}, we have exploited the fact that $\fS$ and $\fU$ are perfect twins because every pair $(v,w), (w,v)$ in the arc set of a symmetric digraph corresponds uniquely to the edge $\{v,w\}$ of an undirected graph and vice versa. In this way, we transferred all concepts and results about directed graphs to undirected graphs, as long as they were compatible with symmetry. Because the class $\fTa \subset \fD$ refers to antisymmetry, we had to replace it by the class $\fC_o \subset \fU$ defined as
\begin{equation*}
\fC_o \; \equiv \; \mysetdescr{ G \in \fU }{ G^* \mytext{\em  does not contain a cycle of odd length} }.
\end{equation*}
Denoting with $\H_u$ and $\S_u$ the sets of homomorphisms and of strict homomorphisms of undirected graphs, Theorem \ref{theo_GschemeOnStrict} became
\begin{theorem}[\cite{aCampo_toappear_0}, Theorem 4] \label{theo_GschemeOnStrict_UG}
Let $R \in \fU$ and $\fU' \subseteq \fU$. Then, for all $S \in \fU$, the equivalence
\begin{align*}
R & \strG S \; \; \mytext{with respect to } \fU'\\ \Leftrightarrow \quad \# \S_u(G,R) & \leq \# \S_u(G,S) \; \; \mytext{for all} \; G \in \fU',
\end{align*}
and the implication
\begin{align*}
\# \S_u(G,R) & \leq \# \S_u(G,S) \; \; \mytext{for all} \; G \in \fU' \\
\Rightarrow \quad \# \H_u(G,R) & \leq \# \H_u(G,S) \; \; \mytext{for all} \; G \in \fU'
\end{align*}
hold if
\begin{align*}
R & \in \fU' = \fU, \\
\mytext{or} \quad R & \in \fC_o \subseteq \fU' \subseteq \fU.
\end{align*}
\end{theorem}

Also for the present article, most of the concepts and results can be directly transferred to undirected graphs just by replacing ``digraph'' by ``undirected graph'', ``$\fD'$'' and ``$\fD$'' by ``$\fU'$ and ``$\fU$'', and ``$\H$'' and ``$\S$ ''by ``$\H_u$'' and ``$\S_u$''. In particular, S-schemes, simple S-schemes, induced S-schemes, and Condition \ref{cond_univ_aexiv} can be defined for undirected graphs in this way, and all results of Section \ref{subsec_simple_induced} can stereotypically be translated into results about undirected graphs. This includes in particular the main results of this section: the characterization of induced S-schemes provided in the Theorems \ref{theo_induced_exists} and \ref{theo_eta_RS}.

Also in the Sections \ref{sec_EVSystems} and \ref{sec_rearr}, the first step is to replace symbols, but now also the concepts change slightly. The definition of the EV-system for undirected graphs has to be modified in an obvious manner. Using the symbol $A(G)$ also for the edge set of an undirected graph $G$, we define
\begin{definition} \label{def_EVsys_alt_UG}
Let $R$ be an undirected graph and $\fU' \subseteq \fU$. We define
\begin{align*}
\E_o(R) & \equiv \mysetdescr{ ( v, D ) }{ v \in V(R), D \subseteq  N_R(v) }.
\end{align*}
For $\fa \in \E_o(R)$, we refer to the two components of $\fa$ by $\fa_1$ and $\fa_2$, and we define
\begin{align*}
\phi_R : \E(R) & \rightarrow R, \\
\fa & \mapsto \fa_1.
\end{align*}
Furthermore, for every $G \in \fU', \xi \in \S_u(G,R)$, we define the mapping
\begin{align*}
\myaTgxi{R}{G}{\xi} : V(G) & \rightarrow \E_o(R), \\
v & \mapsto \left( \xiv, \xi[ N_G(v) ] \right).
\end{align*}
The {\em EV-system $\E(R)$ of $R$ with respect to $\fU'$} is the undirected graph with vertex set $\E_o(R)$ and edge set $A( \E(R) )$ defined by
\begin{align*}
& \{ \fa, \fb \} \in A( \E(R) ) \\
\equiv \quad & \exists \; G \in \fU', \xi \in \S_u(G,R), \{v, w\} \in A(G) \mytext{: } \fa = \aRgxiv, \; \fb = \aRgxiw.
\end{align*}
\end{definition}
Again, $\alpha$ turns out to be a simple S-scheme from $R$ to $\E(R)$, and Lemma \ref{lemma_fafb_simpleProp} becomes
\begin{lemma} \label{lemma_fafb_simpleProp_UG}
For all $\fa, \fb \in \E_o(R)$,
\begin{align*}
\fa \fb \in A( \E(R ) ) \quad & \Rightarrow \quad \fa_1 \fb_1 \in A(R), \\
\fa \fb \in A( \E(R)^* ) \quad & \Rightarrow \quad 
\fa_1 \fb_1 \in A( R^* ), \fa_1 \in \fb_2, \fb_1 \in \fa_2, \\
\fa_1 = \fb_1 \quad & \Rightarrow \quad \fa = \fb.  \\
\mytext{Furthermore,} \quad R \in \fC_o \quad & \Rightarrow \quad \E(R) \in \fC_o.
\end{align*}
\end{lemma}
(For the proof of $\E(R) \in \fC_o$, start with a walk $\fc^0, \ldots , \fc^I$ of odd length in $\E(R)^*$, proceed as in the original, and observe, that in the case of $\fc^0 = \fc^I$, the sequence $\fc^0_1, \ldots , \fc^I_1$ is a cycle of odd length in $R^*$.)

In the following Lemma \ref{lemma_phi_def} and Lemma \ref{lemma_phi_fa}, we just have to skip everything indexed with $3$. Now we define the objects $X(\fa) \in \fC_o$:
\begin{definition} \label{def_Xfa_UG}
For every $m \in \myN_0$, we define the undirected graph $X_m \in \fC_o$ by
\begin{align*}
V( X_m ) & \; \equiv \; D \cup \{ p \}, \\
A( X_m ) & \; \equiv \; \mysetdescr{ \{ p, d \} }{ d \in D },
\end{align*}
where $D$ is a set with $\# D = m$ and $p \notin D$.

Furthermore, for $\fa \in \E_o(R)$, we define
\begin{equation*}
X(\fa) \; \equiv \; X_{\# \fa_2}.
\end{equation*}
$\iota(\fa) : V(X(\fa)) \rightarrow V(R)$ is a mapping sending $p$ to $\fa_1$ and $D$ bijectively to $\fa_2$.
\end{definition}
From now on, we have to proceed more carefully. The reason is the undirected graph $X_1$ defined up to isomorphism by
\begin{align*}
V(X_1) & \; \equiv \; \{ p, d \} \mytext{with } p \not= d, \\
A(X_1) & \; \equiv \; \{ \{ p, d \} \}.
\end{align*}
$\Aut(X_1)$ contains {\em two} automorphisms $\id_{V(X_1)}$ and $\chi$, the latter one interchanging $p$ and $d$. The vertex $p$ is not a fixed point of $\chi$, and we have to introduce a case discrimination in rewriting Corollary \ref{coro_props_Xfa}:
\begin{corollary} \label{coro_props_Xfa_UG}
Let $R$ and $\fU'$ as in the choices in Theorem \ref{theo_GschemeOnStrict_UG}. For every $\fa \in \E_o(R)$, $\pi \in \Aut(X(\fa))$, we have
\begin{align*}
\myaxiv{\iota(\fa)}{p} 
& \; = \; \fa, \\
\myaxiv{\iota(\fa) \circ \pi}{v} 
& \; = \;
\myaiav{\fa}{\pi(v)} \quad \mytext{for all } v \in V(X(\fa)), \pi \in \Aut( X(\fa) ).
\end{align*}
 Furthermore,
\begin{equation*}
\# \Aut( X(\fa ) ) \; = \; \# \mysetdescr{ \iota(\fa) \circ \pi }{ \pi \in \Aut( X(\fa) ) } \; = 
\begin{cases}
\; 2 & \mytext{if } X(\fa) \simeq X_1, \\
\# \fa_2 ! & \mytext{if } X(\fa) \not\simeq X_1.
\end{cases}
\end{equation*}
For $X(\fa) \not\simeq X_1$, the vertex $p$ is a fixed point of every $\pi \in \Aut(X(\fa)$, and $\pi$ maps $D$ bijectively to $D$.
\end{corollary}

To point out the difference: Definition \ref{def_Xfa_UG} assigns the undirected graph $X(\fa) = X_1$ to $\fa = ( p, \{d\}) \in \E_o(X_1)$. However, if we regard $X_1$ as symmetric digraph $Y_1$ with $V(Y_1) = \{ p, d \}$, $A(Y_1) = \{ (p,d), (d,p) \}$, then Definition \ref{def_Xfa} assigns to $\fb = (p,\{d\},\{d\}) \in \E_o(Y_1)$ the {\em antisymmetric} digraph $X(\fb) \simeq X_1^1 \in \fTa$ with {\em three} vertices and $\# \Aut( X_1^1) = 1$, cf.\ Figure \ref{figure_Xij}(b).

Proposition \ref{prop_fDStr_fTa} remains nearly unchanged:

\begin{proposition} \label{prop_fDStr_fTa_UG}
Let $\fC_o \subseteq \fU' \subseteq \fU$ and $R \in \fU$. For all $\fa, \fb \in \E_o(R)$,
\begin{align*}
\fa \fb \in A( \E(R)^* ) \quad & \Leftrightarrow \quad 
\fa_1 \in \fb_2 \mytext{ and } \fb_1 \in \fa_2, \\
\fa \fa \in A( \E(R) ) \quad & \Leftrightarrow \quad \fa_1 \fa_1 \in A(R). 
\end{align*}
ERD and AID hold for $\fU' = \fU$, and for $R \in \fC_o \subseteq \fU'$.
\end{proposition}

Now we come to the counterpart of Theorem \ref{theo_eta_RS_inv}, which we get by just replacing symbols:
\begin{theorem} \label{theo_eta_RS_inv_UG}
Let $R \in \fU' = \fU$ or $R \in \fC_o \subseteq \fU' \subseteq \fU$, and let $\E(R)$ and $\E(S)$ be the EV-Systems  of $R$ and $S \in \fU$. If $\rho$ is a strong S-scheme from $R$ to $S$ fulfilling \eqref{eq_imagebased_gl}, then $\eps \equiv \myaTgxi{S}{\E(R)}{\rho( \phi_R ) }$ fulfills Condition \ref{cond_univ_aexiv}, induces $\rho$, and is one-to-one.
\end{theorem}

However, in the proof, we have to take into account the special status of $X_1$. Corollary \ref{coro_rhophi_DpU} becomes
\begin{corollary} \label{coro_rhophi_DpU_UG}
Let $\fa \in \E_o(R)$. Then, for every $\pi \in \Aut(X(\fa))$,
\begin{align*}
\rho(\iota(\fa) \circ \pi)(v) 
& \; = \;
\rho(\iota(\fa)(\pi(v))
\quad \mytext{for all } v \in V(X(\fa)).
\end{align*}
In the case of $X(\fa) \not\simeq X_1$, we have
\begin{align*}
\begin{split}
\rho(\iota(\fa) \circ \pi)(p) & \; = \; \eps(\fa)_1, \\
\rho(\iota(\fa) \circ \pi)[D] & \; = \; \eps(\fa)_2.
\end{split}
\end{align*}
and in the case of $X(\fa) \simeq X_1$, we have $\eps( \fa ) = ( a, \{ b \} )$ with $ab \in A(S^*)$ and
\begin{align*}
a & \; = \; \myrhoxiv{\iota(\fa)}{p} \; = \; \myrhoxiv{ \iota(\fa) \circ \chi}{ d }, \\
b & \; = \; \myrhoxiv{\iota(\fa)}{d} \; = \; \myrhoxiv{ \iota(\fa) \circ \chi}{ p }.
\end{align*}
\end{corollary}

Lemma \ref{lemma_rhophi_isom_1} remains unchanged:
\begin{lemma} \label{lemma_rhophi_isom_1_UG}
Let $\fa \in \E_o(R)$. If $\rho$ is strong, then, for every $\pi \in \Aut(X(\fa))$, the mapping $\rho( \iota(\fa) \circ \pi)$ is one-to-one on $D$. In particular, $ X( \fa ) \simeq X( \eps( \fa ) )$.
\end{lemma}
(For $X(\fa) \not\simeq X_1$, run the proof as in the original using Corollary \ref{coro_rhophi_DpU_UG}; for $X(\fa) \simeq X_1$, the statement about $D$ is trivial, and $ X( \fa ) \simeq X( \eps( \fa ) )$ holds due to the addendum in Corollary \ref{coro_rhophi_DpU_UG}.)

Finally, in the proof of Lemma \ref{lemma_eps_oneone}, we have to introduce a case disrimination again:

\begin{lemma} \label{lemma_eps_oneone_UG}
If $\rho$ is strong, then $\eps$ is one-to-one.
\end{lemma} 
\BP Let $\fa, \fb \in \E(R)$ with $\eps( \fa ) = \eps( \fb )$. According to Lemma \ref{lemma_rhophi_isom_1_UG}, we have $X(\fa) \simeq X(\eps(\fa)) = X(\eps(\fb)) \simeq X(\fb)$. $X(\fa)$ and $X(\fb)$ are thus isomorphic, and due to $X(\fa), X(\fb) \in \fU'$, we have $X(\fa) = X(\fb)$. Let $G \equiv X(\fa)$ and
\begin{align*}
\I(\fa) & \; \equiv \; \mysetdescr{ \iota(\fa) \circ \pi }{\pi \in \Aut(G)}, \\
\I(\fb) & \; \equiv \; \mysetdescr{ \iota(\fb) \circ \pi }{ \pi \in \Aut(G) }.
\end{align*}
For $X(\fa) \not\simeq X_1$, proceed as in the original. Now assume $X(\fa) \simeq X_1$. According to Corollary \ref{coro_props_Xfa_UG}, the cardinality of $\Aut( G )$, $\I(\fa)$, and $\I(\fb)$ is $2$. $\J \equiv \I(\fa) \cup \I(\fb)$ is a subset of $\S(G, R)$ with $\# \J = 4$ in the case of $\{ \fa_1 \} \cup \fa_2 \not= \{ \fb_1 \} \cup \fb_2$ and $\# \J = 2$ in the case of $\{ \fa_1 \} \cup \fa_2 = \{ \fb_1 \} \cup \fb_2$.

$\rho_G[ \J ]$ is a subset of $\S_u(G,S)$, and due to $\eps(\fa) = \eps(\fb)$ and the last two equations in Corollary \ref{coro_rhophi_DpU_UG}, all elements of $\rho_G[ \J ]$ map $V(G) = \{ p, d \}$ bijectively to the same two-element-subset $\{ v, w \}$ of $V(S)$. We conclude $\# \rho_G[ \J ] \leq 2$, and because $\rho$ is strong, we have $\# \J = \# \rho_G[ \J ]$, hence $\{ \fa_1 \} \cup \fa_2 = \{ \fb_1 \} \cup \fb_2$.

Let $\fa = ( v, \{ w \} )$. If $\fb = ( w, \{ v \} )$, then $\fa \fb \in A(\E(R)^*)$ according to Proposition \ref{prop_fDStr_fTa_UG} in contradiction to $\eps( \fa ) = \eps( \fb )$. Therefore, $\fa = \fb$.

\EP
In Lemma \ref{lemma_aexiv_axiv}, everything indexed with $3$ has to be skiped, and the replacement proposition becomes
\begin{proposition} \label{prop_repl_Cond1_UG}
Let $\eps : \E(R) \rightarrow \E(S)$ be a strict homomorphism between the EV-systems of $R$ and $S$. Assume that for $\fa \in \E_o(R)$
\begin{align*}
\begin{split}
\# \eps( \fa )_2 & \; \leq \; \# \fa_2
\end{split}
\end{align*}
and
\begin{align*}
\forall \; \fb, \fc \in N_{\E(R)}( \fa ) \mytext{: } \eps( \fb )_1 = \eps( \fc )_1 & \; \Rightarrow \; \fb_1 = \fc_1.
\end{align*}
Then $\aSexiv = \eps( \aRxiv )$ for all $G \in \fU'$, $\xi \in \S_u(G,R)$, $v \in V(G)$ with $\aRxiv = \fa$. In particular, $\eps$ fulfills Condition \ref{cond_univ_aexiv} if these conditions hold for all $\fa \in \E_o(R)$.
\end{proposition}

We have already seen in \cite{aCampo_toappear_0}, that also the rearrangement method recalled in Section \ref{subsec_rearr_induced} can be rewritten for undirected graphs. For $R \in \fU$ with the properties described in Definition \ref{def_leq_s}, define $S \in \fU$ by
\begin{align*}
V(S) & \equiv \; Z, \\
A(S) & \equiv \; \left( A(R) \setminus A_{M,X} \right) \; \cup \; A_b, \\
\mytext{where} \quad \quad A_{M,X} & \equiv \; \mysetdescr{ e \in A(R) }{ e \cap M \not= \emptyset \mytext{ and } e \cap X \not= \emptyset} \\
\mytext{and} \quad \quad \quad \; A_b & \equiv \;
\mysetdescr{ ( e \cap M ) \cup \beta[ e \cap X ] }{ e \in A_{M,X} }.
\end{align*}
Then the counterpart of Theorem \ref{theo_rhoxi} delivers $R \strG S$ if we replace \eqref{bed_nbh} by
\begin{align*}
\forall x \in X \mytext{:} & \; N_R(x) \setminus M \subseteq N_R( \beta(x) ).
\end{align*}
With these modifications, all results and formulas in Section \ref{subsec_rearr_induced} remain valid if we just make the usual rewritings and skip everything indexed with $3$.

\end{document}